# EDGE-REINFORCED RANDOM WALK ON A LADDER

By Franz Merkl and Silke W. W. Rolles[1]

*University of Munich and University of Bielefeld*

We prove that the edge-reinforced random walk on the ladder $\mathbb{Z} \times \{1,2\}$ with initial weights $a > 3/4$ is recurrent. The proof uses a known representation of the edge-reinforced random walk on a finite piece of the ladder as a random walk in a random environment. This environment is given by a marginal of a multicomponent Gibbsian process. A transfer operator technique and entropy estimates from statistical mechanics are used to analyze this Gibbsian process. Furthermore, we prove spatially exponentially fast decreasing bounds for normalized local times of the edge-reinforced random walk on a finite piece of the ladder, uniformly in the size of the finite piece.

**1. Introduction.** The oldest reinforced random walk is the edge-reinforced random walk introduced by Diaconis in 1987. The process can be defined on any locally finite nondirected graph as follows: Every edge is given a weight which changes in time. Initially, all edges are given weight 1, say. In each step, the reinforced random walker jumps to a nearest-neighbor vertex with probability proportional to the weight of the traversed edge. Each time an edge is traversed, its weight is increased by 1.

In the late 1980s, Diaconis asked whether the edge-reinforced random walk on $\mathbb{Z}^d$, $d \geq 1$, is recurrent or transient. This problem is still open for all dimensions $d \geq 2$. In one dimension, the edge-reinforced random walk is recurrent. To see this, one can use that on $\mathbb{Z}$, the edge-reinforced random walk has the same distribution as a random walk in an independent random environment. Pemantle [10] proved a phase transition in the recurrence/transience behavior for the edge-reinforced random walk on an infinite binary tree. He used a representation as a random walk in an i.i.d. environment. This method completely fails for graphs with cycles.

Received October 2003; revised July 2004.
[1]Supported in part by NSF Grant DMS-00-71766.

In this article, we prove that the edge-reinforced random walk on $\mathbb{Z} \times \{1, 2\}$ is recurrent. The problem is much more subtle than for acyclic graphs. We use the fact that the edge-reinforced random walk on a finite ladder has the same distribution as a random walk in an environment given by random *time-independent* edge weights. These edge weights are stochastically dependent in a complicated way. The representation as a random walk in a random environment follows from a de Finetti theorem for reversible Markov chains, proven in [12], which refines a result of Diaconis and Freedman [3] in the special case of reversible chains. The distribution of the environment is given by a joint density which was discovered by Coppersmith and Diaconis (see [2]). Our proof relies on a generalization of their statement, proven by Keane and Rolles [7].

It was suggested independently by Diaconis and Keane that the representation as a random walk in a random environment may be useful to prove recurrence for edge-reinforced random walks. Keane conjectured that, in addition, the use of transfer operators may help to prove recurrence on the ladder. The present article pursues such an approach for the first time successfully.

The edge-reinforced random walk studied in this article behaves very differently from the directed-edge-reinforced random walk where every undirected edge is replaced by two directed edges which both get their own weight. The directed-edge-reinforced random walk has the same distribution as a random walk in an *independent* environment (see [8]). This representation is used in [8] to show recurrence on $\mathbb{Z} \times G$ for any finite graph $G$.

A simpler model is the once-reinforced random walk where the weight of an edge is increased by a fixed parameter $\delta > 0$ the first time it is traversed. From the second traversal on, the weight of an edge does not change. This random walk is recurrent on $\mathbb{Z}^1$ (e.g., [1]); however, its recurrence/transience behavior on $\mathbb{Z}^d$ is not known for $d \geq 2$. Even for ladders $\mathbb{Z} \times \{1, 2, \ldots, d\}$ the problem is subtle: For $\delta \in (0, 1/(d-2))$, recurrence was proved by Sellke [14]. Vervoort [16] extended the result for very large $\delta$. For intermediate values of $\delta$ the problem seems to be open. Durrett, Kesten and Limic [6] showed that the once-reinforced random walk on regular trees is transient for all $\delta > 0$.

For the integer line, Davis [1] proved a recurrence/transience dichotomy for a general class of reinforced random walks, including the edge-reinforced and once-reinforced random walk on $\mathbb{Z}$. Very strong localization was shown for random walks with superlinear edge-reinforcement by Limic [9]. Vertex-reinforced random walk localizes as well. This was proved by Pemantle and Volkov [11, 17] and by Tarrès [15].

1.1. *Results.* In this article we consider an edge-reinforced random walk on the ladder $\mathbb{Z} \times \{1, 2\}$. The edges are undirected. They are assigned time-dependent random weights, with all initial edge weights equal to some constant $a > 0$. In each step, the random walker jumps to a nearest-neighbor vertex with probability proportional to the weight of the traversed edge. Whenever the random walk crosses an edge, its weight is increased by 1.

Formally, the edge-reinforced random walk on a locally finite graph $G = (V, E)$ is defined as follows: Let $X_t : V^{\mathbb{N}_0} \to V$ denote the canonical projection on the $t$th coordinate; $X_t$ is interpreted as the random location of the random walker at time $t$. We identify an edge with the set of its endpoints. For $t \in \mathbb{N}_0$, we define $w_t(e) : V^{\mathbb{N}_0} \to \mathbb{R}_+$, the weight of edge $e$ at time $t$, recursively as follows:

$$w_0(e) := a \quad \text{for all } e \in E, \tag{1.1}$$

$$w_{t+1}(e) := \begin{cases} w_t(e) + 1, & \text{for } e = \{X_t, X_{t+1}\} \in E, \\ w_t(e), & \text{for } e \in E \setminus \{\{X_t, X_{t+1}\}\}. \end{cases} \tag{1.2}$$

Let $P^G_{v_0,a}$ denote the distribution of the edge-reinforced random walk on $G$ starting in $v_0$ with all initial edge weights equal to $a$. The distribution $P^G_{v_0,a}$ is a probability measure on $V^{\mathbb{N}_0}$, specified by the following requirements:

$$X_0 = v_0 \quad P^G_{v_0,a}\text{-a.s.}, \tag{1.3}$$

$$\begin{aligned}P^G_{v_0,a}[X_{t+1} = v | X_i, i = 0, 1, \ldots, t] \\ = \begin{cases} \dfrac{w_t(\{X_t, v\})}{\sum_{\{e \in E \,:\, X_t \in e\}} w_t(e)}, & \text{if } \{X_t, v\} \in E, \\ 0, & \text{otherwise.} \end{cases}\end{aligned} \tag{1.4}$$

The ladder is the graph $G = (V, E)$ with vertex set $V := \mathbb{Z} \times \{1, 2\}$ and edge set $E := \{\{u, v\} : u, v \in V \text{ with } \|u - v\|_1 = 1\}$, where $\|\cdot\|_1$ denotes the 1-norm. We say that the reinforced random walk is *recurrent* if almost all paths visit all vertices infinitely often. Our main result reads as follows:

THEOREM 1.1. *For all $a > 3/4$, the edge-reinforced random walk on $\mathbb{Z} \times \{1, 2\}$ with all initial weights equal to $a$ is recurrent.*

The theorem includes the most interesting case $a = 1$. However, we do not expect the bound $3/4$ to be optimal.

Fix $a > 0$ and $n \in \mathbb{N}$. Let $G^{(n)} = (V^{(n)}, E^{(n)})$ with $V^{(n)} := \{0, 1, 2, \ldots, n\} \times \{1, 2\}$ denote the finite ladder of length $n$. We consider an edge-reinforced random walk on $G^{(n)}$ with all initial edge weights equal to $a$. For the vertices, we introduce the following notation (see Figure 1):

$$\underline{i} := (i, 1) \quad \text{and} \quad \overline{i} := (i, 2) \quad \text{for } i \in \mathbb{Z}. \tag{1.5}$$

We abbreviate $P^{(n)} := P^{G^{(n)}}_{\underline{0},a}$. Let $k_t(e)$ denote the (random) number of times the reinforced random walker traverses the (undirected) edge $e$ up to time $t$. We prove:

4  F. MERKL AND S. W. W. ROLLES

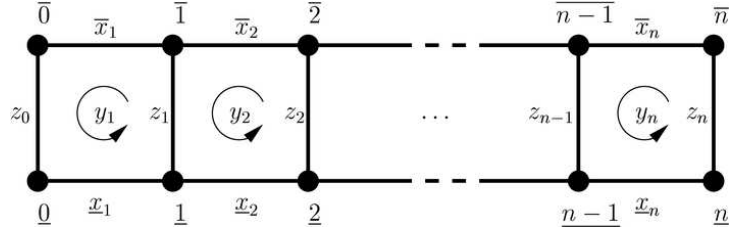

Fig. 1. *The finite ladder.*

THEOREM 1.2. *For all $a > 3/4$ and all $n \in \mathbb{N}$ the following holds*: *If $e \in E^{(n)}$ is an edge on level $i \in \{1, 2, \ldots, n\}$ of the ladder, that is, $\underline{i} \in e$ or $\overline{i} \in e$, then we have*

$$(1.6) \qquad P^{(n)}\left[\lim_{t \to \infty} \frac{k_t(e)}{k_t(\{\underline{0}, \overline{0}\})} \leq e^{-c_1 i}\right] \geq 1 - c_2 e^{-c_3 i}$$

*with constants $c_1, c_2, c_3 > 0$ depending only on $a$.*

The edge-reinforced random walk on a finite graph has the same distribution as a random walk in a *time-independent* random environment. A proof can be found in [12]; see also [4]. For a precise formulation of the result, we refer to Theorem 2.2 below. The time-independent random environment on $G^{(n)}$ is given by *time-independent* random weights $x = (x_e)_{e \in E^{(n)}} \in \mathbb{R}_+^{E^{(n)}}$ with respect to some probability measure $\tilde{\mathbb{Q}}^{(n)}$.

Let us explain how time-independent weights $x = (x_e)_{e \in E^{(n)}} \in \mathbb{R}_+^{E^{(n)}}$ induce a random walk on $G^{(n)}$: In each step, the random walker jumps to a neighboring vertex with probability proportional to the weight $x_e$ of the connecting edge $e$. This random walk is a reversible Markov chain.

The distribution of the edge weights $x = (x_e)_{e \in E^{(n)}}$ constituting the random environment equals the distribution of the limit $(\lim_{t \to \infty} k_t(e)/t)_{e \in E^{(n)}}$, which is $P^{(n)}$-a.s. strictly positive in all components. Therefore, Theorem 1.2 states the following: If we normalize the weights such that the edge $\{\underline{0}, \overline{0}\}$ has weight 1, then the weights of the edges in $E^{(n)}$ decay exponentially in the level of the ladder, uniformly in $n$, with probability close to 1. We use Theorem 1.2 to estimate the escape probability for the reinforced random walker on $G^{(n)}$.

1.2. *Overview of the proofs.* The representation of edge-reinforced random walks as a random walk in a random environment is essential for the whole paper. This representation is described in [12]; we review it in Section 2.1. In order to see the Gibbsian structure behind the random environment $x$, additional auxiliary variables are introduced: These are a random spanning tree $T$ and Gaussian random variables $y$, indexed by a basis of the first homology space of the ladder.

As our third main result, Theorem 2.3 bounds the tails of ratios between components of $x$ and components of $y$. A scaling transformation is used to get rid of the nonlocal constraint that the random environment $x$ is supported on a high-dimensional simplex. The spanning tree $T$ is described by a list of local discrete variables, subject to local matching rules. This is explained in Section 2.2.

The joint probability law of $x$, $y$ and $T$ is Gibbsian. In order to derive estimates for pieces of the corresponding Hamiltonian, we need to transform $x$, $y$ and $T$ to new local variables, described in Section 2.3. Roughly speaking, these new variables consist of logarithms of ratios between neighboring components of $x$ and $y$ and of the local description of the spanning tree $T$. The total Hamiltonian is written in these new coordinates and split into local pieces. The transformation is technically involved.

Lower bounds for the local Hamiltonians are derived in Section 3. The proof of these bounds requires, among other things, the solution of some linear optimization problems. Here, the precise choice of the splitting to local pieces is essential; a "naive" choice would not suffice.

Given the Gibbsian description of the random environment with appropriately bounded local Hamiltonians, we use a transfer operator technique to analyze the random environment. The transfer operator is introduced in Section 4.1; it turns out to be a Hilbert–Schmidt operator.

In Section 4.2 we introduce an auxiliary deformation of the original Gibbs measure which has a reflection symmetry, at least in an asymptotic sense. In the language of statistical mechanics, we apply additional "external forces" at the endpoints of the Gibbsian chain. We show that the Gibbsian chain is deformed by these forces roughly proportional to its length. At the heart of the argument lies an upper bound for the free energy of the deformed Gibbs measure, which is obtained by the variational principle for free energies. In Section 4.3 the deformed Gibbs measure serves to prove exponentially decreasing bounds for the random environment.

Finally, recurrence of the edge-reinforced random walk is proven in Section 5. The exponential bounds for the random environment and the connection between random walks and electric networks are the key ingredients to prove recurrence on the half-sided ladder $\mathbb{N}_0 \times \{1, 2\}$. Symmetry and gluing arguments imply recurrence on the two-sided ladder.

**2. Edge-reinforced random walk on the finite ladder.** Fix $a > 0$ and $n \in \mathbb{N}$. In this section we consider the edge-reinforced random walk on the finite ladder $G^{(n)}$ of length $n$. All initial edge weights are equal to $a$, and the random walker starts in $\overline{0}$.

We could equally well study the walk on $\{-n, \ldots, n\} \times \{1, 2\}$ started at $\overline{0}$ in the middle, rather than the walk on $\{0, \ldots, n\} \times \{1, 2\}$ started at the left boundary. However, the version presented here is simpler, since it avoids distinguishing three cases "to the left of the starting point," "at the starting point" and "to the right of the starting point" for inner points of the ladder in many estimates presented below.

If there is no risk of confusion, we omit sometimes the dependence on $n$ in the notation.

### 2.1. A random walk in a random environment.
We need some notation: Let

$$\Delta^{(n)} := \left\{ x = (x_e)_{e \in E^{(n)}} \in (0,1]^{E^{(n)}} : \sum_{e \in E^{(n)}} x_e = 1 \right\} \tag{2.1}$$

denote the simplex. We define $\alpha_t := (k_t(e)/t)_{e \in E^{(n)}}$. Clearly, $\alpha_t \in \Delta^{(n)}$.

For $v \in V$, we denote by $x_v$ the sum of all $x_e$ with $e$ incident to $v$:

$$x_v := \sum_{\{e \in E : v \in e\}} x_e. \tag{2.2}$$

Furthermore, we introduce the following abbreviations:

$$\underline{x}_i := x_{\{i-1,\underline{i}\}}, \qquad \overline{x}_i := x_{\{\overline{i-1},\overline{i}\}}, \qquad z_i := x_{\{\underline{i},\overline{i}\}}. \tag{2.3}$$

Figure 1 illustrates this definition.

For $1 \leq i \leq n$, let $c_i$ denote the oriented cycle $\underline{i-1}$, $\underline{i}$, $\overline{i}$, $\overline{i-1}$, $\underline{i-1}$. Clearly, $c_1, \ldots, c_n$ constitute a basis of the first homology space $H_1(G^{(n)})$ (cycle space) of $G^{(n)}$; note that the first Betti number $\dim H_1(G^{(n)})$ of $G^{(n)}$ equals $n$. Let $\mathbb{R}_+ := (0, \infty)$. For $x \in \mathbb{R}_+^{E^{(n)}}$, we define the matrix $A^{(n)}(x) = (A^{(n)}_{i,j}(x))_{1 \leq i,j \leq n}$ by

$$A^{(n)}_{i,i}(x) := \sum_{e \in c_i} \frac{1}{x_e}, \qquad A^{(n)}_{i,j}(x) := \sum_{e \in c_i \cap c_j} \pm \frac{1}{x_e} \qquad \text{for } i \neq j, \tag{2.4}$$

where the signs in the last sum are chosen to be $+1$ or $-1$ depending on whether the edge $e$ has in $c_i$ and $c_j$ the same orientation or not. Explicitly, this means

$$\begin{aligned}
A^{(n)}_{i,i}(x) &= \frac{1}{z_{i-1}} + \frac{1}{\underline{x}_i} + \frac{1}{\overline{x}_i} + \frac{1}{z_i}, \\
A^{(n)}_{i,i+1}(x) &= -\frac{1}{z_i} = A^{(n)}_{i+1,i}(x), \\
A^{(n)}_{i,j}(x) &= 0 \qquad \text{for } |i-j| \geq 2.
\end{aligned} \tag{2.5}$$

Let $\mathcal{T}^{(n)}$ denote the set of all spanning trees of $G^{(n)}$. Let $y^t$ denote the transpose of a vector $y$, and let $E(T)$ denote the edge set of a tree $T$. For

$x = (x_e)_{e \in E^{(n)}} \in \mathbb{R}_+^{E^{(n)}}$, $y = (y_1, y_2, \ldots, y_n) \in \mathbb{R}^n$ and $T \in \mathcal{T}^{(n)}$, we define

$$\Phi^{(n)}(x, y, T) = \frac{\prod_{i=1}^n [\underline{x}_i^{a-3/2} \overline{x}_i^{a-3/2}] \prod_{i=0}^n z_i^{a-3/2} \prod_{e \in E(T)} x_e}{x_{\underline{0}}^{a+1/2} x_{\overline{0}}^a \prod_{i=1}^{n-1} [x_{\underline{i}}^{(3a+1)/2} x_{\overline{i}}^{(3a+1)/2}] x_{\underline{n}}^{a+1/2} x_{\overline{n}}^{a+1/2}}$$

(2.6)

$$\times \exp\left[-\frac{1}{2} y A^{(n)}(x) y^t\right].$$

Let $\sigma$ denote the Lebesgue measure on $\Delta^{(n)}$, normalized so that $\sigma(\Delta^{(n)}) = 1$. We set

(2.7) $$\tilde{z}^{(n)} := \int_{\Delta^{(n)}} \int_{\mathbb{R}^n} \sum_{T \in \mathcal{T}^{(n)}} \Phi^{(n)}(x, y, T) \, dy \, \sigma(dx).$$

The normalizing constant $\tilde{z}^{(n)}$ is given explicitly in Theorem 1 of [7]; in particular it is finite.

THEOREM 2.1 ([7], Theorem 1). *The sequence $(\alpha_t)_{t \in \mathbb{N}}$ converges almost surely. The distribution of the limit is absolutely continuous with respect to the surface measure $\sigma$ on $\Delta^{(n)}$ with density given by*

(2.8) $$\phi^{(n)}(x) = \frac{1}{\tilde{z}^{(n)}} \int_{\mathbb{R}^n} \sum_{T \in \mathcal{T}^{(n)}} \Phi^{(n)}(x, y, T) \, dy.$$

There is also a probabilistic interpretation of the arguments $y$ in (2.6) in terms of winding numbers of the reinforced random walk paths; for details see [7].

The edge-reinforced random walk on every finite graph has the same distribution as a random walk in a random environment where the environment is given by weights on the edges. We state the result only for $G^{(n)}$:

THEOREM 2.2 ([12], Theorem 3.1). *For any path $(v_0, v_1, \ldots, v_t)$ in $G^{(n)}$, the following holds:*

(2.9) $$P^{(n)}[X_i = v_i \text{ for } 0 \leq i \leq t] = \int_{\Delta^{(n)}} \prod_{i=1}^t \frac{x_{\{v_{i-1}, v_i\}}}{x_{v_{i-1}}} \phi^{(n)}(x) \, \sigma(dx);$$

*here $x := (x_e)_{e \in E^{(n)}}$. Hence, $P^{(n)}$ equals the distribution of the random walk in a random environment on $G^{(n)}$ starting in $\overline{0}$ with environment given by random edge weights chosen according to $\phi^{(n)} \, d\sigma$.*

We define

(2.10)
$$\tilde{\Lambda}^{(n)} := \Delta^{(n)} \times \mathbb{R}^n \times \mathcal{T}^{(n)},$$
$$d\tilde{\mathbb{Q}}^{(n)} := [\tilde{z}^{(n)}]^{-1} \Phi^{(n)}(x, y, T) \, \sigma(dx) \, dy \, dT,$$

where $dT$ denotes the counting measure on $\mathcal{T}^{(n)}$. The marginal of the distribution $\tilde{\mathbb{Q}}^{(n)}$ with respect to the components $x$ equals the distribution of the random environment as a measure on $\Delta^{(n)}$.

The following theorem bounds the tails of ratios between these random variables and the "winding number" random variable $y_i^2$.

THEOREM 2.3. *Let $a > 3/4$. With respect to $\tilde{\mathbb{Q}}^{(n)}$, all the random variables*

$$(2.11) \qquad \ln \frac{\underline{x}_i}{y_i^2}, \qquad \ln \frac{\overline{x}_i}{y_i^2}, \qquad \ln \frac{z_i}{y_i^2} \quad and \quad \ln \left| \frac{y_{i+1}}{y_i} \right|$$

*have exponential tails, uniformly in $i$ and $n$. In other words, there exist constants $c_4(a) > 0$ and $c_5(a) > 0$ such that for all $n \in \mathbb{N}$, $i < n$ and $M > 0$, one has*

$$(2.12) \qquad \tilde{\mathbb{Q}}^{(n)}[|\Upsilon_i| \geq M] \leq c_4 e^{-c_5 M},$$

*where $\Upsilon_i$ denotes any of the four random variables in (2.11).*

For $x = (x_e)_{e \in E^{(n)}} \in \Delta^{(n)}$, all weights $(cx_e)_{e \in E^{(n)}}$ with $c > 0$ induce the same reversible Markov chain. In Theorem 2.2, the edge weights $x$ are normalized in such a way that $\sum_e x_e = 1$. For our purposes, it is more convenient to set one weight, namely $z_0$, equal to 1; recall that we used this normalization in Theorem 1.2. The change of normalization is made precise in Lemma 2.6. We set

$$(2.13) \qquad \Lambda^{(n)} := \mathbb{R}_+^{E^{(n)}} \times \mathbb{R}^n \times \mathcal{T}^{(n)}.$$

The following scaling property of $\Phi^{(n)}$ will be important. We omit its elementary proof.

LEMMA 2.4. *Let $c > 0$ be a real number. For all $(x, y, T) \in \Lambda^{(n)}$, the following holds:*

$$(2.14) \qquad \Phi^{(n)}(cx, c^{1/2}y, T) = c^{-(7/2)n-1} \Phi^{(n)}(x, y, T).$$

For $x \in \mathbb{R}_+^{E^{(n)}}$, we set $s(x) := \sum_{e \in E^{(n)}} x_e$. We define

$$(2.15) \qquad S: \Lambda^{(n)} \to \tilde{\Lambda}^{(n)}, \qquad (x, y, T) \mapsto (s(x)^{-1}x, s(x)^{-1/2}y, T).$$

We write $|A|$ for the cardinality of a set $A$. Let $\delta_1$ denote the Dirac measure in 1.

DEFINITION 2.5. *Let $\lambda(dx) := \delta_1(dz_0) \times \prod_{i=1}^n d\underline{x}_i \, d\overline{x}_i \, dz_i$. We define*

$$(2.16) \qquad d\mathbb{Q}^{(n)} := \frac{(3n)!}{\tilde{z}(n)} \Phi^{(n)}(x, y, T) \lambda(dx) \, dy \, dT.$$

LEMMA 2.6. $\mathbb{Q}^{(n)}$ is the image measure of $\tilde{\mathbb{Q}}^{(n)}$ under $S$; that is, for any measurable function $f:\tilde{\Lambda}^{(n)} \to \mathbb{R}_+$, the following holds:

$$\text{(2.17)} \qquad \int_{\tilde{\Lambda}^{(n)}} f\, d\tilde{\mathbb{Q}}^{(n)} = \int_{\Lambda^{(n)}} f \circ S\, d\mathbb{Q}^{(n)}.$$

In particular, for any path $(v_0, v_1, \ldots, v_t)$ in $G^{(n)}$, we have

$$\text{(2.18)} \quad P^{(n)}[X_i = v_i \text{ for } 0 \leq i \leq t] = \int_{\Lambda^{(n)}} \prod_{i=1}^{t} \frac{x_{\{v_{i-1}, v_i\}}}{x_{v_{i-1}}}\, d\mathbb{Q}^{(n)}(x, y, T).$$

PROOF. Let $f: \tilde{\Lambda}^{(n)} \to \mathbb{R}_+$ be measurable. Then, introducing an auxiliary integration over $t$ and using the definition (2.10) of $\tilde{\mathbb{Q}}^{(n)}$,

$$\begin{aligned}
\int_{\tilde{\Lambda}^{(n)}} f\, d\tilde{\mathbb{Q}}^{(n)} &= \int_{\tilde{\Lambda}^{(n)}} \int_0^\infty e^{-t} f(x', y', T)\, dt\, d\tilde{\mathbb{Q}}^{(n)}(x', y', T) \\
\text{(2.19)} \qquad &= \frac{1}{\tilde{z}^{(n)}} \int_{\tilde{\Lambda}^{(n)}} \int_0^\infty t^{-7n/2} e^{-t} f(x', y', T) \\
&\qquad \times \Phi^{(n)}(x', y', T)\, t^{7n/2}\, dt\, \sigma(dx')\, dy'\, dT.
\end{aligned}$$

Consider the transformation $\Lambda^{(n)} \to (0, \infty) \times \tilde{\Lambda}^{(n)}$, defined by $(x, y, T) \mapsto (t, x', y', T) = (s(x), S(x, y, T))$. Then, the measure $t^{7n/2}\, dt\, \sigma(dx')\, dy'\, dT$ on the right-hand side of (2.19) is the image measure of the measure $(3n)!\, dx\, dy\, dT$ under this transformation. To see this, note first that the projection to $|E^{(n)}| - 1 = 3n$ components $(x'_e)_{e \in E^{(n)} \setminus \{e_0\}}$ suffices to parametrize any point $(x'_e)_{e \in E^{(n)}} \in \Delta^{(n)}$, due to the constraint $\sum_{e \in E^{(n)}} x'_e = 1$; here we abbreviate $e_0 := \{\underline{0}, \overline{0}\}$. Second, the image of the probability measure $\sigma$ on the simplex $\Delta^{(n)}$ under this projection $(x'_e)_{e \in E^{(n)}} \mapsto (x'_e)_{e \in E^{(n)} \setminus \{e_0\}}$ equals $(3n)!$ times the Lebesgue measure on $\{(x'_e)_e \in \mathbb{R}_+^{E^{(n)} \setminus \{e_0\}} : \sum_{e \neq e_0} x_e < 1\}$; the normalizing factor $(3n)! = (|E^{(n)}| - 1)!$ arises since the last set has the volume $1/(|E^{(n)}| - 1)!$. Third, the inverse transformation, with the redundant component $x'_{e_0}$ and the discrete variable $T$ dropped, is given by $(t, (x'_e)_{e \neq e_0}, y') \mapsto ((tx'_e)_{e \neq e_0}, t(1 - \sum_{e \neq e_0} x'_e), \sqrt{t}y') = ((x_e)_{e \neq e_0}, x_{e_0}, y)$. Its Jacobi determinant is given by $t^{|E^{(n)}| - 1}(\sqrt{t})^n = t^{7n/2}$. Hence, by the transformation formula, the right-hand side of (2.19) equals

$$\begin{aligned}
&\frac{(3n)!}{\tilde{z}^{(n)}} \int_{\Lambda^{(n)}} s(x)^{-7n/2} e^{-s(x)} f(S(x, y, T)) \Phi^{(n)}(S(x, y, T))\, dx\, dy\, dT \\
\text{(2.20)} \qquad &= \frac{(3n)!}{\tilde{z}^{(n)}} \int_{\Lambda^{(n)}} s(x) e^{-s(x)} f(S(x, y, T)) \Phi^{(n)}(x, y, T)\, dx\, dy\, dT;
\end{aligned}$$

for the last equality we used Lemma 2.4 with $c = s(x)^{-1}$. Next, we substitute $\tilde{x} = z_0^{-1} x$, $\tilde{y} = z_0^{-1/2} y$. Note that the $e_0$-component of $\tilde{x}$ equals 1. Consider the transformation $(x, y, T) \mapsto (z_0, \tilde{x}, \tilde{y}, T)$. If we drop the redundant component $\tilde{x}_{e_0} = 1$ and the discrete variable $T$, then the resulting transformation $(x, y) \mapsto (z_0, (\tilde{x}_e)_{e \neq e_0}, \tilde{y})$ has the Jacobi determinant $z_0^{-(|E^{(n)}|-1)} z_0^{-n/2} = z_0^{-7n/2}$. Thus, the image measure of $dx\, dy\, dT$ with respect to the transformation $(x, y, T) \mapsto (z_0, \tilde{x}, \tilde{y}, T)$ equals $z_0^{7n/2}\, dz_0\, \lambda(d\tilde{x})\, d\tilde{y}\, dT$; note that the Dirac measure in $\lambda(d\tilde{x})$ arises from $\tilde{x}_{e_0} = 1$. Using $S(x, y, T) = S(\tilde{x}, \tilde{y}, T)$ and once more the transformation formula, the expression (2.20) becomes

$$\text{(2.21)} \quad \frac{(3n)!}{\tilde{z}^{(n)}} \int_{\Lambda^{(n)}} \int_0^\infty z_0 s(\tilde{x}) e^{-z_0 s(\tilde{x})} f(S(\tilde{x}, \tilde{y}, T)) \\ \times \Phi^{(n)}(z_0 \tilde{x}, z_0^{1/2} \tilde{y}, T) z_0^{7n/2}\, dz_0\, \lambda(d\tilde{x})\, d\tilde{y}\, dT.$$

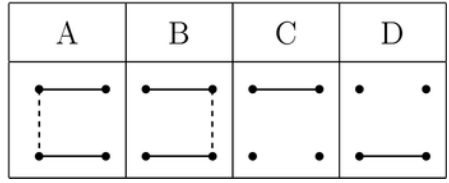

Fig. 2. *The possible states of the tree variables.*

By Lemma 2.4, $\Phi^{(n)}(z_0 \tilde{x}, z_0^{1/2} \tilde{y}, T) z_0^{7n/2+1} = \Phi^{(n)}(\tilde{x}, \tilde{y}, T)$. Hence, integrating over $z_0$ yields for the expression (2.21):

$$\text{(2.22)} \quad \frac{(3n)!}{\tilde{z}^{(n)}} \int_{\Lambda^{(n)}} f(S(\tilde{x}, \tilde{y}, T)) \Phi^{(n)}(\tilde{x}, \tilde{y}, T)\, \lambda(d\tilde{x})\, d\tilde{y}\, dT \\ = \int_{\Lambda^{(n)}} f \circ S\, d\mathbb{Q}^{(n)}.$$

This completes the proof of (2.17). The statement (2.18) follows immediately from Theorem 2.2 because the quotients $x_{\{v_{i-1}, v_i\}}/x_{v_{i-1}}$ do not change under the transformation $S$. □

2.2. *Building spanning trees.* In this section we give a local description of the spanning trees of $G^{(n)}$. In order to describe a spanning tree $T$, we specify for each of the $n$ cycles one of the states $A$, $B$, $C$ or $D$ shown in Figure 2.

Let $\mathbb{T}^{(n)} := \{(T_i)_{1 \leq i \leq n} \in \{A, B, C, D\}^n : (T_i, T_{i+1}) \neq (A, B)$ for all $i = 1, 2, \ldots, n-1\}$, and recall that $\mathcal{T}^{(n)}$ denotes the set of all spanning trees of $G^{(n)}$. We define

$$\text{(2.23)} \quad \Psi_{\text{tree}} : \mathbb{T}^{(n)} \to \mathcal{T}^{(n)}, \qquad \Psi_{\text{tree}}((T_i)_{1 \leq i \leq n}) := T;$$

$T_i$ describes which of the horizontal edges $\{\underline{i-1}, \underline{i}\}$, $\{\overline{i-1}, \overline{i}\}$ is contained in the spanning tree $T$ according to Figure 2. Let us define which rungs are included in $T$. The left rung $\{\underline{0},\overline{0}\}$ is included for $T_1 \in \{A,C,D\}$. Figure 3 tells us for $i=1,\ldots,n-1$ whether the rung $\{\underline{i},\overline{i}\}$ is contained in the tree $T$. Finally, the right rung $\{\underline{n},\overline{n}\}$ is included for $T_n \in \{B,C,D\}$.

Indeed, $\Psi_{\text{tree}}((T_i)_{1\leq i\leq n})$ yields a spanning tree of $G^{(n)}$.

LEMMA 2.7. *The map* $\Psi_{\text{tree}} \colon \mathbb{T}^{(n)} \to \mathcal{T}^{(n)}$ *is a bijection.*

PROOF. Let $T \in \mathcal{T}^{(n)}$. For each of the cycles $c_i$ ($1 \leq i \leq n$) of $G^{(n)}$, there is at least one edge which is not contained in $T$. If $\{\underline{i-1},\underline{i}\} \notin E(T)$, we set $T_i := C$. If $\{\overline{i-1},\overline{i}\} \notin E(T)$, we set $T_i := D$. Otherwise, both edges are contained in $T$, and $T$ is connected by a vertical edge somewhere either on the left (and we set $T_i := A$) or on the right (and $T_i := B$). Clearly, $(T_i, T_{i+1}) \neq (A,B)$ for all $i$. One sees that $\Psi_{\text{tree}}$ maps the constructed sequence $(T_i)_{1\leq i\leq n}$ to $T$. Hence, $\Psi_{\text{tree}}$ is onto, and it is not hard to see that it is one-to-one as well. □

2.3. *A Gibbsian representation of the random environment.* In this section we represent $\mathbb{Q}^{(n)}$ as the image measure of a Gibbsian probability measure under a suitable transformation. This representation is essential for the analysis of the random environment. The transformed local variables are called $\underline{X}_i$, $\overline{X}_i$, $\sigma_i$, $T_i$, $Z_i$ and $\Gamma_i$. In the terminology of statistical mechanics, one may view them as abstract local "spin" variables. The new variables $\underline{X}_i, \overline{X}_i, Z_i, \Gamma_i$ consist essentially of logarithms of ratios between neighboring old variables $x_e$, $y_i^2$, while $\sigma_i$ and $T_i$ are discrete. First, we describe the space of these transformed local variables.

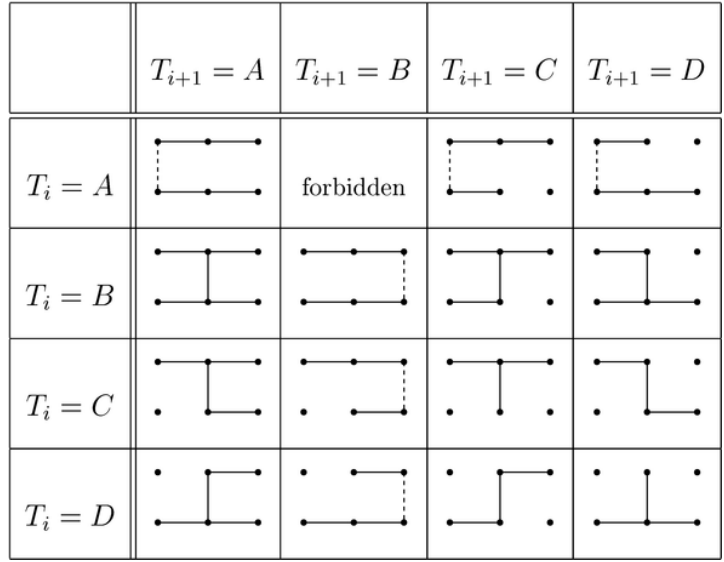

FIG. 3. *Possible transitions of the tree variables.*

DEFINITION 2.8. We define $\Omega^{(n)} := \Omega_{\text{left}} \times \Omega_{\text{cycle}}^n \times \Omega_{\text{rung}}^{n-1} \times \Omega_{\text{right}}$, where

$$\begin{aligned}\Omega_{\text{left}} &= \mathbb{R}, & \Omega_{\text{cycle}} &= \mathbb{R}^2 \times \{\pm 1\} \times \{A, B, C, D\}, \\ \Omega_{\text{rung}} &= \mathbb{R}^2, & \Omega_{\text{right}} &= \mathbb{R}.\end{aligned} \quad (2.24)$$

We endow $\Omega^{(n)}$ with the reference measure $d\omega$, defined as the Lebesgue measure on $\mathbb{R}^{4n}$ times the counting measure. We denote the canonical projections on $\Omega^{(n)}$ by

$$(2.25) \quad \omega_{\text{left}}, \quad (\omega_{\text{cycle},i})_{1 \leq i \leq n}, \quad (\omega_{\text{rung},i})_{1 \leq i \leq n-1}, \quad \omega_{\text{right}},$$

where $\omega_{\text{left}} = Z_0$, $\omega_{\text{cycle},i} = (\underline{X}_i, \overline{X}_i, \sigma_i, T_i)$, $\omega_{\text{rung},i} = (Z_i, \Gamma_i)$, $\omega_{\text{right}} = Z_n$; hence a generic element on $\Omega^{(n)}$ is

$$(2.26) \quad \omega = (Z_0, (\underline{X}_i, \overline{X}_i, \sigma_i, T_i)_{1 \leq i \leq n}, (Z_i, \Gamma_i)_{1 \leq i \leq n-1}, Z_n).$$

Furthermore, we set

$$\begin{aligned}(2.27) \quad \tilde{\Omega}^{(n)} := \{\omega \in \Omega^{(n)} &: (T_i(\omega), T_{i+1}(\omega)) \neq (A, B) \\ &\text{for all } i = 1, \ldots, n-1\}.\end{aligned}$$

Intuitively speaking, the components of $\Omega^{(n)}$ are associated to parts of the ladder: $\Omega_{\text{left}}$ and $\Omega_{\text{right}}$ belong to the left and right rung, respectively; the $i$th $\Omega_{\text{cycle}}$ component belongs to the $i$th cycle in the ladder, and the $i$th $\Omega_{\text{rung}}$ component is associated to the step from the $i$th to the $(i+1)$st cycle.

Using an intuitive statistical mechanics picture, one has a chain of "compound spins," consisting of "inner spin variables" $\underline{X}_i$, $\overline{X}_i$, $Z_i$, $\sigma_i$ and $T_i$, while $\Gamma_i$ measures the "separation" between neighboring compound spins.

This is clarified by the following definition; it describes the transformation to the new local variables. In order to simplify the notation we first define auxiliary variables $U_i$, $W_i$ and $Y_i$.

DEFINITION 2.9. We introduce the abbreviations

$$(2.28) \quad U_i := \tfrac{1}{2}[\underline{X}_i + \overline{X}_i], \qquad 1 \leq i \leq n,$$

$$(2.29) \quad W_i := \Gamma_i + U_{i+1} - U_i, \qquad 1 \leq i \leq n-1,$$

$$(2.30) \quad Y_i := -Z_0 - \sum_{j=1}^{i-1} W_j, \qquad 1 \leq i \leq n.$$

Recall the abbreviations (2.3). We set

$$(2.31) \quad \Lambda_1^{(n)} := \{((x_e)_e, (y_i)_i, T) \in \Lambda^{(n)} : z_0 = 1,\ y_i \neq 0 \text{ for all } i\},$$

and we define $\Psi : \tilde{\Omega}^{(n)} \to \Lambda_1^{(n)}$, $\omega \mapsto ((x_e)_{e \in E^{(n)}}, (y_i)_{1 \leq i \leq n}, T)$, as follows:

(2.32)
$$\underline{x}_i = e^{\underline{X}_i(\omega)+Y_i(\omega)}, \qquad \overline{x}_i = e^{\overline{X}_i(\omega)+Y_i(\omega)},$$
$$y_i = \sigma_i e^{Y_i(\omega)/2}, \qquad 1 \leq i \leq n;$$

(2.33)
$$z_0 = 1, \qquad z_i = e^{Z_i(\omega)+(Y_i(\omega)+Y_{i+1}(\omega))/2},$$
$$z_n = e^{Z_n(\omega)+Y_n(\omega)}, \qquad 1 \leq i \leq n-1;$$

(2.34) $\qquad T = \Psi_{\text{tree}}(T_1(\omega), \ldots, T_n(\omega)).$

For our analysis below, it is essential to consider as our new (capital) variables logarithms of quotients roughly of the type $x_e/x_{e'}$ or $x_e/y_i^2$, with $e$ close to $e'$ and close to level $i$, since these variables turn out to have uniformly exponential tails. There is some arbitrariness in the choice of the precise form of the change of variables. Our choice was optimized so that the estimate for the local Hamiltonians in Proposition 3.2, is valid for sufficiently small values of $a$.

Speaking very roughly and using the statistical mechanics picture again, $Y_i$ may be viewed as the center of the $i$th compound spin, while $\underline{X}_i$, $\overline{X}_i$ and $Z_i$ can be considered as distances of the constituents of the $i$th compound spin to the center. The constituents are bound together by strong internal forces. In the intuitive picture, the locations $\underline{X}_i + Y_i$, $\overline{X}_i + Y_i$ and $Z_i + \frac{1}{2}(Y_i + Y_{i+1})$ of the constituents are logarithms of the original edge weights $x_e$. Multiplicative scaling of these weights corresponds to shifting the whole spin chain.

We note that $\Lambda_1^{(n)}$ has full $\mathbb{Q}^{(n)}$-measure, that is, $\mathbb{Q}^{(n)}[\Lambda_1^{(n)}] = 1$, and that $\Psi : \tilde{\Omega}^{(n)} \to \Lambda_1^{(n)}$ is a bijection. Its inverse $\Psi^{-1} : \Lambda_1^{(n)} \to \tilde{\Omega}^{(n)}$ is expressed by the following equations:

(2.35)
$$\underline{X}_i = \ln \frac{\underline{x}_i}{y_i^2}, \qquad \overline{X}_i = \ln \frac{\overline{x}_i}{y_i^2}, \qquad \sigma_i = \operatorname{sgn} y_i,$$
$$T_i = [\Psi_{\text{tree}}^{-1}(T)]_i, \qquad 1 \leq i \leq n,$$

(2.36) $\qquad Z_0 = \ln \dfrac{z_0}{y_1^2}, \qquad Z_i = \ln \dfrac{z_i}{|y_i y_{i+1}|}, \qquad Z_n = \ln \dfrac{z_n}{y_n^2},$

(2.37) $\qquad \Gamma_i = \dfrac{1}{2}\left[\ln \dfrac{\underline{x}_i}{\underline{x}_{i+1}} + \ln \dfrac{\overline{x}_i}{\overline{x}_{i+1}}\right], \qquad 1 \leq i \leq n-1.$

The following relations are also useful:

(2.38)
$$Y_i = 2\ln|y_i|$$
$$= \tfrac{1}{2}(\underline{X}_1 + \overline{X}_1 - \underline{X}_i - \overline{X}_i) - Z_0 - \sum_{j=1}^{i-1} \Gamma_j, \qquad 1 \leq i \leq n,$$

$$(2.39) \quad W_i = Y_i - Y_{i+1} = \ln \frac{y_i^2}{y_{i+1}^2},$$

$$(2.40) \quad \Gamma_i = [U_i + Y_i] - [U_{i+1} + Y_{i+1}], \qquad 1 \le i \le n-1.$$

The reason to take $y_i$ squared is its scaling behavior in (2.14).

The image measure of $\mathbb{Q}^{(n)}$ with respect to the transformation $\Psi^{-1}: \Lambda_1^{(n)} \to \tilde{\Omega}^{(n)} \subset \Omega^{(n)}$ turns out to be a Gibbs measure; we show this in Lemma 2.13. But first, we introduce the relevant local Hamiltonians for this Gibbs measure. Since we deform the Gibbs measure later, we introduce a "deformation parameter" $\eta$ already at this point. In the language of statistical mechanics, one may view $\eta$ as an external force, coupling to the "separation" $\Gamma_i$ between compound spins. Without deformation, $\eta$ takes the value $1/4$.

In analogy to (2.28) and (2.29), we use the following abbreviations:

$$(2.41) \quad U = \tfrac{1}{2}(\underline{X} + \overline{X}), \qquad U' = \tfrac{1}{2}(\underline{X}' + \overline{X}'), \qquad W = \Gamma + U' - U.$$

DEFINITION 2.10. For $\eta \in \mathbb{R}$, we define $H_{\text{middle},a,\eta}: \Omega_{\text{cycle}} \times \Omega_{\text{rung}} \times \Omega_{\text{cycle}} \to \mathbb{R} \cup \{\infty\}$,

$$H_{\text{middle},a,\eta}(\underline{X}, \overline{X}, \sigma, T | Z, \Gamma | \underline{X}', \overline{X}', \sigma', T')$$
$$(2.42) \quad := H_{\ln,a} + H_{\text{linear},a} + H_{\text{tree}} + H_{\text{expI}} + H_{\text{expII}} + H_{\text{constraint}} - \eta \Gamma$$
$$= H_{\text{middle},a,0} - \eta \Gamma,$$

where

$$(2.43) \quad H_{\ln,a} := \frac{3a+1}{2} \{\ln[e^{\underline{X}+W/2} + e^{\underline{X}'-W/2} + e^Z]$$
$$+ \ln[e^{\overline{X}+W/2} + e^{\overline{X}'-W/2} + e^Z]\},$$

$$(2.44) \quad H_{\text{linear},a} := -(a+\tfrac{1}{2})[U + U' + Z],$$

$$H_{\text{tree}} := \tfrac{1}{2}[\mathbb{1}_{\{T=C\}} \underline{X} + \mathbb{1}_{\{T=D\}} \overline{X} + \mathbb{1}_{\{T'=C\}} \underline{X}' + \mathbb{1}_{\{T'=D\}} \overline{X}']$$
$$(2.45) \quad + [\mathbb{1}_{\{T'=B\}} + \mathbb{1}_{\{T=A\}}]Z + \tfrac{1}{2}[\mathbb{1}_{\{T'=B\}} - \mathbb{1}_{\{T=A\}}]W$$
$$- \tfrac{1}{2}[\mathbb{1}_{\{T=A\}} - \mathbb{1}_{\{T=B\}}]U + \tfrac{1}{2}[\mathbb{1}_{\{T'=A\}} - \mathbb{1}_{\{T'=B\}}]U',$$

$$(2.46) \quad H_{\text{expI}} := \tfrac{1}{4}[e^{-\underline{X}} + e^{-\overline{X}} + e^{-\underline{X}'} + e^{-\overline{X}'}],$$

$$(2.47) \quad H_{\text{expII}} := \tfrac{1}{2}[\sigma e^{W/4} - \sigma' e^{-W/4}]^2 e^{-Z},$$

$$(2.48) \; H_{\text{constraint}} := +\infty \cdot \mathbb{1}_{\{T=A, T'=B\}}.$$

The piece $H_{\text{constraint}}$ encodes the constraint $(T_i, T_{i+1}) \ne (A, B)$ for the tree variables; intuitively speaking, the energy required to violate this constraint is infinite.

We define the reflection $\leftrightarrow$ by $A^\leftrightarrow = B$, $B^\leftrightarrow = A$, $C^\leftrightarrow = C$ and $D^\leftrightarrow = D$. For $\omega = (\underline{X}, \overline{X}, \sigma, T) \in \Omega_{\text{cycle}}$, we define $\omega^\leftrightarrow := (\underline{X}, \overline{X}, \sigma, T^\leftrightarrow)$. We observe the following reflection symmetry property:

$$(2.49) \qquad H_{\text{middle},a,\eta}(\omega|Z,\Gamma|\omega') = H_{\text{middle},a,-\eta}(\omega'^\leftrightarrow|Z,-\Gamma|\omega^\leftrightarrow).$$

The two ends of the ladder need some extra consideration:

DEFINITION 2.11. We define $H_{\text{left},a} : \Omega_{\text{left}} \times \Omega_{\text{cycle}} \to \mathbb{R}$ by

$$(2.50) \quad H_{\text{left},a}(Z|\underline{X}, \overline{X}, \sigma, T) = H_{\text{left},a} := H_{\text{left},\ln,a} + H_{\text{left,tree}} + H_{\exp} + \frac{U}{4},$$

where

$$(2.51) \quad H_{\text{left},\ln,a} := a\ln[e^{\overline{X}} + e^Z] + (a+\tfrac{1}{2})\{\ln[e^{\underline{X}} + e^Z] - [U+Z]\},$$

$$(2.52) \quad \begin{aligned} H_{\text{left,tree}} &:= \tfrac{1}{2}[\mathbb{1}_{\{T=C\}}\underline{X} + \mathbb{1}_{\{T=D\}}\overline{X}] \\ &\quad + \mathbb{1}_{\{T=B\}}Z + \tfrac{1}{2}[\mathbb{1}_{\{T=A\}} - \mathbb{1}_{\{T=B\}}]U, \end{aligned}$$

$$(2.53) \qquad H_{\exp} := \tfrac{1}{4}[e^{-\underline{X}} + e^{-\overline{X}}] + \tfrac{1}{2}e^{-Z}.$$

We define $H_{\text{right},a} : \Omega_{\text{cycle}} \times \Omega_{\text{right}} \to \mathbb{R}$ by

$$(2.54) \quad \begin{aligned} H_{\text{right},a}(\underline{X}, \overline{X}, \sigma, T|Z) \\ = H_{\text{right},a} := H_{\text{right},\ln,a} + H_{\text{right,tree}} + H_{\exp} - \frac{U}{4} \end{aligned}$$

with

$$(2.55) \quad H_{\text{right},\ln,a} := (a+\tfrac{1}{2})\{\ln[e^{\underline{X}} + e^Z] + \ln[e^{\overline{X}} + e^Z] - [U+Z]\},$$

$$(2.56) \quad \begin{aligned} H_{\text{right,tree}} &:= \tfrac{1}{2}[\mathbb{1}_{\{T=C\}}\underline{X} + \mathbb{1}_{\{T=D\}}\overline{X}] \\ &\quad + \mathbb{1}_{\{T=A\}}Z - \tfrac{1}{2}[\mathbb{1}_{\{T=A\}} - \mathbb{1}_{\{T=B\}}]U. \end{aligned}$$

The total Hamiltonian is defined to be the sum of all local pieces; here the deformation parameter $\eta$ equals $1/4$:

DEFINITION 2.12. We define $H^{(n)} : \Omega^{(n)} \to \mathbb{R} \cup \{\infty\}$ by

$$(2.57) \quad \begin{aligned} H^{(n)}(\omega) &:= H_{\text{left},a}(\omega_{\text{left}}|\omega_{\text{cycle},1}) \\ &\quad + \sum_{i=1}^{n-1} H_{\text{middle},a,1/4}(\omega_{\text{cycle},i}|\omega_{\text{rung},i}|\omega_{\text{cycle},i+1}) \\ &\quad + H_{\text{right}}(\omega_{\text{cycle},n}|\omega_{\text{right}}). \end{aligned}$$

The following lemma shows that the Gibbs measure defined by the total Hamiltonian $H^{(n)}$ indeed describes the random environment for the reinforced random walk, transformed to the local spin variables.

Recall Definition 2.9 of the bijection $\Psi$ and the definition of the reference measure $d\omega$ on $\Omega^{(n)}$ in Definition 2.8.

LEMMA 2.13. *Let $\mathbb{P}^{(n)}$ denote the image measure of $\mathbb{Q}^{(n)}$ (restricted to the set $\Lambda_1^{(n)}$ of full measure) under the map $\Psi^{-1} : \Lambda_1^{(n)} \to \Omega^{(n)}$, that is, for all measurable functions $f : \Lambda^{(n)} \to \mathbb{R}_+$, we have*

$$\int_{\Lambda^{(n)}} f \, d\mathbb{Q}^{(n)} = \int_{\tilde{\Omega}^{(n)}} f \circ \Psi \, d\mathbb{P}^{(n)}. \tag{2.58}$$

*Then $\mathbb{P}^{(n)}$ has the following representation as a Gibbs measure:*

$$d\mathbb{P}^{(n)} = [Z^{(n)}]^{-1} e^{-H^{(n)}(\omega)} \, d\omega, \tag{2.59}$$

*with $Z^{(n)} := \int_{\Omega^{(n)}} e^{-H^{(n)}(\omega)} \, d\omega$.*

PROOF. The presence of $H_{\text{constraint}}$ in the Hamiltonian guarantees that $e^{-H^{(n)}}$ is indeed supported on $\text{range}(\Psi^{-1}) = \tilde{\Omega}^{(n)} \subset \Omega^{(n)}$, that is, $e^{-H^{(n)}} = 0$ as soon as one of the constraints $(T_i, T_{i+1}) \neq (A, B)$ is violated.

For $f : \Lambda^{(n)} \to \mathbb{R}_+$ as above, we get the following by the transformation formula and the definition (2.16) of $\mathbb{Q}^{(n)}$:

$$\int_{\tilde{\Omega}^{(n)}} f \circ \Psi \, d\mathbb{P}^{(n)} = \int_{\Lambda_1^{(n)}} f \, d\mathbb{Q}^{(n)}$$
$$= \frac{(3n)!}{\tilde{z}^{(n)}} \int_{\tilde{\Omega}^{(n)}} (f \circ \Psi)(\Phi^{(n)} \circ \Psi) \mathcal{J} \, d\omega, \tag{2.60}$$

with an appropriate Jacobi determinant $\mathcal{J}$, specified in (2.70). Recall that the components of $\omega$, described in (2.26), consist of both continuous and discrete parts. The Jacobi determinant $\mathcal{J}$ is taken only with respect to the continuous variables with the constant value $z_0 = 1$ dropped.

Because of (2.60), it suffices to show that

$$-\ln(\Phi^{(n)}(\Psi(\omega))\mathcal{J}(\omega)) = H^{(n)}(\omega) + c_6 \tag{2.61}$$

holds for all $\omega$ in the domain $\tilde{\Omega}^{(n)}$ of $\Psi$ with a constant $c_6(n) \in \mathbb{R}$. We rewrite the density $\Phi^{(n)}$ given in (2.6) as follows:

$$\Phi^{(n)}(x, y, T)$$
$$= \exp\left[-\frac{1}{2} y A^{(n)}(x) y^t\right] \frac{[\underline{x}_1 \overline{x}_1 z_0^2]^{a/2+1/4}}{x_{\underline{0}}^{a+1/2} x_{\overline{0}}^{a}} \frac{[\underline{x}_n \overline{x}_n z_n^2]^{a/2+1/4}}{[x_{\underline{n}} x_{\overline{n}}]^{a+1/2}} \tag{2.62}$$
$$\times \prod_{i=1}^{n-1} \frac{[\underline{x}_i \overline{x}_i \underline{x}_{i+1} \overline{x}_{i+1} z_i^2]^{a/2+1/4}}{[x_{\underline{i}} x_{\overline{i}}]^{(3a+1)/2}} \prod_{i=1}^{n} [\underline{x}_i \overline{x}_i]^{-2} \prod_{i=0}^{n} z_i^{-2} \prod_{e \in E(T)} x_e.$$

Using the explicit form (2.5) of $A^{(n)}$, we express the argument of the first exponential term as follows:

$$\frac{1}{2}yA^{(n)}(x)y^t = \frac{1}{2}\left\{\frac{y_1^2}{z_0} + \frac{y_n^2}{z_n} + \sum_{i=1}^{n}\left[\frac{y_i^2}{\underline{x}_i} + \frac{y_i^2}{\overline{x}_i}\right] + \sum_{i=1}^{n-1}\frac{(y_i - y_{i+1})^2}{z_i}\right\}$$

$$= \frac{1}{2}\left\{e^{-Z_0} + e^{-Z_n} + \sum_{i=1}^{n}[e^{-\underline{X}_i} + e^{-\overline{X}_i}]\right.$$

(2.63)
$$\left.+ \sum_{i=1}^{n-1}[\sigma_i e^{W_i/4} - \sigma_{i+1}e^{-W_i/4}]^2 e^{-Z_i}\right\}$$

$$= H_{\exp}(\omega_{\text{left}}|\omega_{\text{cycle},1})$$

$$+ \sum_{i=1}^{n-1}(H_{\text{expI}} + H_{\text{expII}})(\omega_{\text{cycle},i}|\omega_{\text{rung},i}|\omega_{\text{cycle},i+1})$$

$$+ H_{\exp}(\omega_{\text{cycle},n}|\omega_{\text{right}}).$$

This is the exponential part of the Hamiltonian $H^{(n)}$. Next, we analyze the part of the density $\Phi^{(n)}$ which depends on the parameter $a$. For $1 \leq i \leq n-1$, we have by Definition 2.9 and (2.39)

$$-\ln\frac{[\underline{x}_i\overline{x}_i\underline{x}_{i+1}\overline{x}_{i+1}z_i^2]^{a/2+1/4}}{[x_{\underline{i}}x_{\overline{i}}]^{(3a+1)/2}}$$

(2.64)
$$= \frac{3a+1}{2}\ln[e^{\underline{X}_i+Y_i} + e^{\underline{X}_{i+1}+Y_{i+1}} + e^{Z_i+(Y_i+Y_{i+1})/2}]$$

$$+ \frac{3a+1}{2}\ln[e^{\overline{X}_i+Y_i} + e^{\overline{X}_{i+1}+Y_{i+1}} + e^{Z_i+(Y_i+Y_{i+1})/2}]$$

$$- \left(a+\frac{1}{2}\right)\left[U_i + U_{i+1} + Z_i + \frac{3}{2}Y_i + \frac{3}{2}Y_{i+1}\right]$$

$$= \frac{3a+1}{2}\ln[e^{\underline{X}_i+W_i/2} + e^{\underline{X}_{i+1}-W_i/2} + e^{Z_i}]$$

$$+ \frac{3a+1}{2}\ln[e^{\overline{X}_i+W_i/2} + e^{\overline{X}_{i+1}-W_i/2} + e^{Z_i}]$$

$$- \left(a+\frac{1}{2}\right)[U_i + U_{i+1} + Z_i] - \frac{1}{4}[Y_i + Y_{i+1}]$$

$$= (H_{\ln,a} + H_{\text{linear},a})(\omega_{\text{cycle},i}|\omega_{\text{rung},i}|\omega_{\text{cycle},i+1}) - \frac{1}{4}[Y_i + Y_{i+1}].$$

Similarly, a calculation using Definition 2.9 and in particular $Y_1 = -Z_0$ yields

$$-\ln\frac{[\underline{x}_1\overline{x}_1z_0^2]^{a/2+1/4}}{x_{\underline{0}}^{a+1/2}x_{\overline{0}}^{a}}$$

$$= \left(a + \frac{1}{2}\right) \ln[e^{\underline{X}_1+Y_1} + 1]$$

(2.65)
$$+ a \ln[e^{\overline{X}_1+Y_1} + 1] - \left(\frac{a}{2} + \frac{1}{4}\right)(\underline{X}_1 + \overline{X}_1 + 2Y_1)$$

$$= a \ln[e^{\overline{X}_1} + e^{Z_0}] + \left(a + \frac{1}{2}\right)\{\ln[e^{\underline{X}_1} + e^{Z_0}] - U_1 - Z_0\} - \frac{Y_1}{2}$$

$$= H_{\text{left},\ln,a}(\omega_{\text{left}}|\omega_{\text{cycle},1}) - \frac{Y_1}{2}$$

and

$$-\ln \frac{[\underline{x}_n \overline{x}_n z_n^2]^{a/2+1/4}}{[x_{\underline{n}} x_{\overline{n}}]^{a+1/2}}$$

$$= \left(a + \frac{1}{2}\right)\{\ln[e^{\underline{X}_n+Y_n} + e^{Z_n+Y_n}] + \ln[e^{\overline{X}_n+Y_n} + e^{Z_n+Y_n}]\}$$

(2.66)
$$- \left(\frac{a}{2} + \frac{1}{4}\right)(\underline{X}_n + \overline{X}_n + 2Z_n + 4Y_n)$$

$$= \left(a + \frac{1}{2}\right)\{\ln[e^{\underline{X}_n} + e^{Z_n}] + \ln[e^{\overline{X}_n} + e^{Z_n}] - U_n - Z_n\}$$

$$= H_{\text{right},\ln,a}(\omega_{\text{cycle},n}|\omega_{\text{right}}).$$

These are the parts of the Hamiltonian which depend on the parameter $a$.

Recall the definition of $\Psi$ from (2.32)–(2.34):

(2.67)
$$\Psi(Z_0, (\underline{X}_i, \overline{X}_i, \sigma_i, T_i)_{1\le i\le n}, (Z_i, \Gamma_i)_{1\le i\le n-1}, Z_n)$$
$$= ((\underline{x}_i, \overline{x}_i)_{1\le i\le n}, z_0, (z_i)_{1\le i\le n-1}, z_n, (y_i)_{1\le i\le n}, T)$$
$$= ((e^{\underline{X}_i+Y_i}, e^{\overline{X}_i+Y_i})_{1\le i\le n}, 1, (e^{Z_i+(Y_i+Y_{i+1})/2})_{1\le i\le n-1},$$
$$e^{Z_n+Y_n}, (\sigma_i e^{Y_i/2})_{1\le i\le n}, \Psi_{\text{tree}}((T_i)_{1\le i\le n}))$$

with

(2.68)
$$Y_i = -Z_0 + \tfrac{1}{2}(\underline{X}_1 + \overline{X}_1 - \underline{X}_i - \overline{X}_i) - \sum_{j=1}^{i-1} \Gamma_j$$

by (2.38). Consider the map $\Psi$ with the discrete arguments $\sigma_i$ and $T_i$ taken fixed and the discrete component $T$ and the constant component $z_0 = 1$ in the image vector dropped. We write this map as a composition of three maps:

$$(Z_0, \underline{X}_1, \overline{X}_1, \ldots, \underline{X}_n, \overline{X}_n, Z_1, \Gamma_1, \ldots, Z_{n-1}, \Gamma_{n-1}, Z_n)$$
$$\stackrel{\pi}{\mapsto} (\underline{X}_1, \overline{X}_1, \ldots, \underline{X}_n, \overline{X}_n, Z_1, \ldots, Z_n, Z_0, \Gamma_1, \ldots, \Gamma_{n-1})$$

(2.69)
$$\stackrel{\chi}{\longmapsto} (\underline{X}_1, \overline{X}_1, \ldots, \underline{X}_n, \overline{X}_n, Z_1, \ldots, Z_n, Y_1, \ldots, Y_n)$$
$$\stackrel{\Xi}{\longmapsto} (\underline{x}_1, \overline{x}_1, \ldots, \underline{x}_n, \overline{x}_n, z_1, \ldots, z_n, y_1, \ldots, y_n).$$

The map $\pi$ is just a permutation; thus, its Jacobi determinant has absolute value 1. Using (2.68), the Jacobi matrix $D\chi$ is a lower triangular matrix having only the entries $+1$ and $-1$ on the diagonal; thus its determinant has also absolute value 1. Finally, the Jacobi matrix $D\Xi$ is an upper triangular matrix; the values on the diagonal are given by $\partial \underline{x}_i/\partial \underline{X}_i = e^{\underline{X}_i + Y_i}$, $\partial \overline{x}_i/\partial \overline{X}_i = e^{\overline{X}_i + Y_i}$, $\partial z_j/\partial Z_j = e^{Z_j + (Y_j + Y_{j+1})/2}$, $\partial z_n/\partial Z_n = e^{Z_n + Y_n}$ and $\partial y_i/\partial Y_i = \frac{1}{2}\sigma_i e^{Y_i/2}$ with $1 \le i \le n$ and $1 \le j \le n-1$.

We observe that the Jacobi determinant $\mathcal{J}$ occurring in (2.60) equals
$$\mathcal{J} = |\det D(\Xi \circ \chi \circ \pi)| = |\det D\Xi|$$

(2.70)
$$= e^{Z_n + Y_n} \prod_{i=1}^n [\tfrac{1}{2} e^{Y_i/2} e^{\underline{X}_i + Y_i} e^{\overline{X}_i + Y_i}] \prod_{i=1}^{n-1} e^{Z_i + (Y_i + Y_{i+1})/2}$$
$$= 2^{-n} \prod_{i=1}^n [e^{Y_i/2} \underline{x}_i \overline{x}_i] \prod_{i=0}^n z_i.$$

Consequently,

(2.71)
$$-\ln\left[\mathcal{J} \prod_{i=1}^n [\underline{x}_i \overline{x}_i]^{-1} \prod_{i=0}^n z_i^{-1}\right] = n \ln 2 - \tfrac{1}{2} \sum_{i=1}^n Y_i.$$

It remains to analyze the part of the density $\Phi^{(n)}$ which depends on the spanning tree $T$. We note that

(2.72)
$$\prod_{i=1}^n [\underline{x}_i \overline{x}_i]^{-1} \prod_{i=0}^n z_i^{-1} \prod_{e \in E(T)} x_e = \prod_{e \notin E(T)} x_e^{-1}.$$

Furthermore, using the definition of the tree variables $A$, $B$, $C$ and $D$ (see Figure 2) and the possible transitions between them (see Figure 3), we obtain

(2.73)
$$-\ln\left[\prod_{e \notin E(T)} x_e^{-1}\right] = \sum_{i=1}^n [\mathbb{1}_{\{T_i = C\}} \ln \underline{x}_i + \mathbb{1}_{\{T_i = D\}} \ln \overline{x}_i]$$
$$+ \sum_{i=1}^{n-1} [\mathbb{1}_{\{T_i = A\}} + \mathbb{1}_{\{T_{i+1} = B\}}] \ln z_i$$
$$+ \mathbb{1}_{\{T_1 = B\}} \ln z_0 + \mathbb{1}_{\{T_n = A\}} \ln z_n$$
$$= \sum_{i=1}^n [\mathbb{1}_{\{T_i = C\}} \underline{X}_i + \mathbb{1}_{\{T_i = D\}} \overline{X}_i]$$

$$+ \sum_{i=1}^{n-1} [\mathbb{1}_{\{T_i=A\}} + \mathbb{1}_{\{T_{i+1}=B\}}]Z_i + \mathbb{1}_{\{T_n=A\}}Z_n$$

$$+ \tfrac{1}{2} \sum_{i=1}^{n-1} [\mathbb{1}_{\{T_{i+1}=B\}} - \mathbb{1}_{\{T_i=A\}}]W_i + R.$$

The remaining term $R$ is specified in (2.74); recall that $W_i = Y_i - Y_{i+1}$ by (2.39). Using the identity $\mathbb{1}_{\{T_i=A\}} + \mathbb{1}_{\{T_i=B\}} + \mathbb{1}_{\{T_i=C\}} + \mathbb{1}_{\{T_i=D\}} = 1$, we see that

$$R = \sum_{i=1}^{n} [\mathbb{1}_{\{T_i=C\}} + \mathbb{1}_{\{T_i=D\}}]Y_i + \sum_{i=1}^{n-1} [\mathbb{1}_{\{T_i=A\}} + \mathbb{1}_{\{T_{i+1}=B\}}]\frac{Y_i + Y_{i+1}}{2}$$

$$+ \mathbb{1}_{\{T_n=A\}}Y_n + \sum_{i=1}^{n-1} [\mathbb{1}_{\{T_{i+1}=B\}} - \mathbb{1}_{\{T_i=A\}}]\frac{Y_{i+1} - Y_i}{2}$$

$$(2.74) \quad = \sum_{i=1}^{n} [1 - \mathbb{1}_{\{T_i=A\}} - \mathbb{1}_{\{T_i=B\}}]Y_i$$

$$+ \sum_{i=1}^{n-1} [\mathbb{1}_{\{T_i=A\}}Y_i + \mathbb{1}_{\{T_{i+1}=B\}}Y_{i+1}] + \mathbb{1}_{\{T_n=A\}}Y_n$$

$$= \sum_{i=1}^{n} Y_i + \mathbb{1}_{\{T_1=B\}}Z_0.$$

For the last equality, we used $Y_1 = -Z_0$ by (2.30). Combining (2.73) and (2.74) with the definitions (2.45), (2.52) and (2.56) of $H_{\text{tree}}$, $H_{\text{left,tree}}$ and $H_{\text{right,tree}}$, we get

$$(2.75) \quad -\ln\left[\prod_{e \notin E(T)} x_e^{-1}\right] = H_{\text{left,tree}}(\omega_{\text{left}}|\omega_{\text{cycle},1})$$

$$+ \sum_{i=1}^{n-1} H_{\text{tree}}(\omega_{\text{cycle},i}|\omega_{\text{rung},i}|\omega_{\text{cycle},i+1})$$

$$+ H_{\text{right,tree}}(\omega_{\text{cycle},n}|\omega_{\text{right}}) + \sum_{i=1}^{n} Y_i.$$

Note that the terms involving $U_i$ in $H_{\text{left,tree}}$, $H_{\text{tree}}$ and $H_{\text{right,tree}}$ form a telescopic sum, which vanishes. The right-hand side in (2.75) is the part of the Hamiltonian depending on the spanning tree $T$.

Next, we sum up (2.63)–(2.66), (2.71) and (2.75). Recall that in the definitions (2.50) of $H_{\text{left},a}$ and (2.54) of $H_{\text{right},a}$, there are additional terms

$\pm U/4$. We get the following result on the set $\tilde{\Omega}^{(n)}$:

$$-\ln((\Phi^{(n)} \circ \Psi)|\mathcal{J}|)$$

$$= H^{(n)} + \frac{1}{4}\sum_{i=1}^{n-1}\Gamma_i - \frac{U_1}{4} + \frac{U_n}{4}$$

(2.76)
$$-\frac{1}{4}\sum_{i=1}^{n-1}(Y_i + Y_{i+1}) - \frac{Y_1}{2} + n\ln 2 - \frac{1}{2}\sum_{i=1}^{n}Y_i + \sum_{i=1}^{n}Y_i$$

$$= H^{(n)} + n\ln 2 + \frac{1}{4}\sum_{i=1}^{n-1}(\Gamma_i - [U_i + Y_i] + [U_{i+1} + Y_{i+1}])$$

$$= H^{(n)} + n\ln 2,$$

where we used (2.40) in the last step. This proves the claim (2.61). □

**3. Bounds for the Hamiltonian.** In this section we bound the local Hamiltonians from below, showing that they increase at least linearly at infinity. Unfortunately, these estimates are technically involved; but fortunately, for the remainder, it suffices to know only the statements of Propositions 3.2 and 3.5 and of Lemma 3.3.

3.1. *Bounding $H_{\text{middle},a,\eta}$.* The key to the linear bound of $H_{\text{middle},a,\eta}$ lies in the solution of a linear optimization problem. It seems to be more convenient to consider a dualized problem; this is why we introduce "dual variables" $\kappa_{T,T'}$ in the following lemma.

LEMMA 3.1. *Let $a = 1/2$. Then there are $\underline{\kappa}_{T,T'} \geq 0$, $\overline{\kappa}_{T,T'} \geq 0$, $\underline{\kappa}'_{T,T'} \geq 0$ and $\overline{\kappa}'_{T,T'} \geq 0$ with*

(3.1) $\qquad\underline{\kappa}_{T,T'} + \overline{\kappa}_{T,T'} = \frac{1}{4}\quad\text{and}\quad \underline{\kappa}'_{T,T'} + \overline{\kappa}'_{T,T'} = \frac{1}{4},$

*such that one has*

(3.2)
$$H_{\ln,1/2} + H_{\text{linear},1/2} + H_{\text{tree}} + H_{\text{constraint}} - \frac{\Gamma}{4}$$
$$\geq \underline{\kappa}_{T,T'}\underline{X} + \overline{\kappa}_{T,T'}\overline{X} + \underline{\kappa}'_{T,T'}\underline{X}' + \overline{\kappa}'_{T,T'}\overline{X}'.$$

PROOF. Due to the presence of $H_{\text{constraint}}$ in the claim (3.2), there is nothing to show in the case $(T,T') = (A,B)$. Thus let us assume $(T,T') \in \{A,B,C,D\}^2 \setminus \{(A,B)\}$. Using $\ln(e^x + e^y + e^z) \geq \max\{x,y,z\}$, we estimate $H_{\ln,1/2}$, defined in (2.43), from below by a convex combination: For all $\underline{\alpha}, \underline{\beta}, \underline{\gamma}, \overline{\alpha}, \overline{\beta}, \overline{\gamma} \geq 0$ with

(3.3) $\qquad\underline{\alpha} + \underline{\beta} + \underline{\gamma} = 1\quad\text{and}\quad \overline{\alpha} + \overline{\beta} + \overline{\gamma} = 1,$

we get

$$H_{\ln,1/2} \geq \frac{5}{4}\left[\max\left\{\underline{X}+\frac{W}{2}, \underline{X}'-\frac{W}{2}, Z\right\} + \max\left\{\overline{X}+\frac{W}{2}, \overline{X}'-\frac{W}{2}, Z\right\}\right]$$

$$(3.4) \quad \geq \frac{5}{4}\left[\underline{\alpha}\left(\underline{X}+\frac{W}{2}\right) + \underline{\beta}\left(\underline{X}'-\frac{W}{2}\right)\right.$$

$$\left. + \underline{\gamma}Z + \overline{\alpha}\left(\overline{X}+\frac{W}{2}\right) + \overline{\beta}\left(\overline{X}'-\frac{W}{2}\right) + \overline{\gamma}Z\right] =: K.$$

It suffices to show that for all 15 possible values $(T,T') \in \{A,B,C,D\}^2 \setminus \{(A,B)\}$ of the tree variables, there exist

$$(3.5) \quad \begin{aligned} &\underline{\alpha} \geq 0, \quad \underline{\beta} \geq 0, \quad \underline{\gamma} \geq 0, \quad \overline{\alpha} \geq 0, \quad \overline{\beta} \geq 0, \quad \overline{\gamma} \geq 0, \\ &\underline{\kappa}_{T,T'} \geq 0, \quad \overline{\kappa}_{T,T'} \geq 0, \quad \underline{\kappa}'_{T,T'} \geq 0, \quad \overline{\kappa}'_{T,T'} \geq 0, \end{aligned}$$

so that (3.1), (3.3) and

$$(3.6) \quad K + H_{\text{linear},1/2} + H_{\text{tree}} - \frac{\Gamma}{4} = \underline{\kappa}_{T,T'}\underline{X} + \overline{\kappa}_{T,T'}\overline{X} + \underline{\kappa}'_{T,T'}\underline{X}' + \overline{\kappa}'_{T,T'}\overline{X}'$$

hold for all real values $\underline{X}$, $\overline{X}$, $\underline{X}'$, $\overline{X}'$, $Z$ and $\Gamma$. Indeed, the claim (3.2) follows from (3.6) and the bound (3.4).

Equation (3.6) is equivalent to the system of equations for the coefficients of $\underline{X}$, $\overline{X}$, $\underline{X}'$, $\overline{X}'$, $Z$ and $\Gamma$ in (3.6). This system together with (3.1), (3.3) and the inequalities (3.5) form a system of linear equations and linear inequalities. We solved this system with a computer for all 15 possibilities for the tree variables $(T,T')$. Here is one possible solution:

$$(3.7) \quad \begin{aligned} \underline{\alpha} := &\tfrac{1}{10}[\mathbb{1}_{\{T'=A\}} - \mathbb{1}_{\{T'=B\}} + \mathbb{1}_{\{T'=C\}} \\ &- \mathbb{1}_{\{T=A,T'=D\}} + \mathbb{1}_{\{T=B,T'=D\}} + \mathbb{1}_{\{T=C,T'=D\}}] \\ &+ \tfrac{1}{5}[3 + \mathbb{1}_{\{T=A\}} - \mathbb{1}_{\{T=B\}} - 2\cdot\mathbb{1}_{\{T=C\}} + \mathbb{1}_{\{T=C,T'=B\}}], \end{aligned}$$

$$(3.8) \quad \begin{aligned} \underline{\beta} := &\tfrac{1}{10}[\mathbb{1}_{\{T'=B\}} - \mathbb{1}_{\{T'=A\}} - 3\cdot\mathbb{1}_{\{T'=C\}} + \mathbb{1}_{\{T=A,T'=D\}}] \\ &+ \tfrac{1}{5}[2 - \mathbb{1}_{\{T=A\}} + \mathbb{1}_{\{T=B\}} - \mathbb{1}_{\{T=B,T'=A\}} + \mathbb{1}_{\{T=C,T'=B\}} \\ &+ \mathbb{1}_{\{T=A,T'=C\}} - \mathbb{1}_{\{T=B,T'=C\}} - \mathbb{1}_{\{T=B,T'=D\}}], \end{aligned}$$

(3.9) $\underline{\gamma} := 1 - \underline{\alpha} - \underline{\beta}$,

(3.10) $\overline{\alpha} := \tfrac{4}{5} + \tfrac{4}{5}\mathbb{1}_{\{T=A\}} - \underline{\alpha}$,

(3.11) $\overline{\beta} := \tfrac{2}{5} + \tfrac{4}{5}\mathbb{1}_{\{T'=B\}} - \underline{\beta}$,

(3.12) $\overline{\gamma} := 1 - \overline{\alpha} - \overline{\beta}$

and

(3.13)
$$\underline{\kappa}_{T,T'} := \tfrac{1}{4}[1 + \mathbb{1}_{\{T=C,T'=B\}} - \mathbb{1}_{\{T'=B\}}$$
$$- \mathbb{1}_{\{T=A,T'=D\}} - \tfrac{1}{2}\mathbb{1}_{\{T=D,T'=D\}}],$$

(3.14)
$$\overline{\kappa}_{T,T'} := \tfrac{1}{4} - \underline{\kappa}_{T,T'},$$

(3.15)
$$\underline{\kappa}'_{T,T'} := \tfrac{1}{4}[1 - \mathbb{1}_{\{T'=B\}} - \mathbb{1}_{\{T=A\}} + \mathbb{1}_{\{T=B,T'=B\}}$$
$$+ \mathbb{1}_{\{T=C,T'=B\}} + \mathbb{1}_{\{T=A,T'=C\}}]$$
$$+ \tfrac{1}{8}[\mathbb{1}_{\{T=A,T'=D\}} - \mathbb{1}_{\{T'=D\}}],$$

(3.16)
$$\overline{\kappa}'_{T,T'} := \tfrac{1}{4} - \underline{\kappa}'_{T,T'}.$$

Indeed, an elementary but tedious calculation shows that (3.6) and the inequalities (3.5) are fulfilled for our choice. We checked all the conditions with the help of a computer algebra system, but they can also be checked by hand for the 15 combinations of the tree variables. This completes the proof. □

The following lemma bounds the local Hamiltonians "in the middle" of the ladder. We split off the exponential part $H_{\text{expII}}$, since this refined bound is needed in a deformation argument in Section 4.2.

PROPOSITION 3.2. *For all $a > 1/2$ and for all $\eta \in [-1/4, 1/4]$ one has*

(3.17)
$$H_{\text{middle},a,\eta} \geq H_{\text{middle},a,\eta} - H_{\text{expII}}$$
$$\geq c_7(a)[|\underline{X}| + |\overline{X}| + |Z| + |\Gamma| + |\underline{X}'| + |\overline{X}'|]$$

*with the constant $c_7(a) = \tfrac{1}{16}\min\{a - \tfrac{1}{2}, 1\} > 0$.*

PROOF. In the estimate (3.19), we use the elementary bound

(3.18)
$$3\ln(e^x + e^y + e^z) - x - y - z \geq 3\max\{x,y,z\} - x - y - z$$
$$\geq \tfrac{1}{2}[|x - y| + |x - z| + |y - z|]$$

twice, for $x = \underline{X} + W/2$, $y = \underline{X}' - W/2$, $z = Z$ and for $x = \overline{X} + W/2$, $y = \overline{X}' - W/2$, $z = Z$. Furthermore, we use an average between the two bounds $|x_1| + |x_2| + |x_3| + |x_4| \geq |x_1 + x_2 + x_3 + x_4|$ and $|x_1| + |x_2| + |x_3| + |x_4| \geq |x_1 - x_3| + |x_2 - x_4|$ in the third step in (3.19), and we abbreviate $\tilde{a} = a - 1/2$. This yields

$$(H_{\ln,a} + H_{\text{linear},a}) - (H_{\ln,1/2} + H_{\text{linear},1/2})$$
$$= \frac{\tilde{a}}{2}[3\ln(e^{\underline{X}+W/2} + e^{\underline{X}'-W/2} + e^Z)$$

$$+ 3\ln(e^{\overline{X}+W/2} + e^{\overline{X}'-W/2} + e^Z) - \underline{X} - \underline{X}' - \overline{X} - \overline{X}' - 2Z]$$

$$\geq \frac{\tilde{a}}{4}\Big[|\underline{X} - \underline{X}' + W| + \Big|\underline{X} + \frac{W}{2} - Z\Big| + \Big|\underline{X}' - \frac{W}{2} - Z\Big|$$

(3.19)
$$+ |\overline{X} - \overline{X}' + W| + \Big|\overline{X} + \frac{W}{2} - Z\Big| + \Big|\overline{X}' - \frac{W}{2} - Z\Big|\Big]$$

$$\geq \frac{\tilde{a}}{4}[|\underline{X} - \underline{X}' + W| + |\overline{X} - \overline{X}' + W|]$$

$$+ \frac{\tilde{a}}{8}[|\underline{X} - \overline{X}| + |\underline{X}' - \overline{X}'| + |\underline{X} + \overline{X} + \underline{X}' + \overline{X}' - 4Z|]$$

$$\geq \frac{\tilde{a}}{4}|\underline{X} + \overline{X} - \underline{X}' - \overline{X}' + 2W|$$

$$+ \frac{\tilde{a}}{8}[|\underline{X} - \overline{X}| + |\underline{X}' - \overline{X}'| + |\underline{X} + \overline{X} + \underline{X}' + \overline{X}' - 4Z|].$$

Lemma 3.1 implies that

(3.20)
$$H_{\ln,1/2} + H_{\text{linear},1/2} + H_{\text{tree}} + H_{\text{constraint}} - \frac{\Gamma}{4}$$
$$\geq \frac{1}{4}[\min\{\underline{X}, \overline{X}\} + \min\{\underline{X}', \overline{X}'\}]$$

holds. We also need the following fact: For all $\varepsilon > 0$, we have

(3.21) $$\min\{x, y\} + e^{-x} + e^{-y} + \varepsilon|x - y| \geq \delta(|x| + |y|),$$

where $\delta = \min\{\varepsilon, 1/2\}$. To prove this fact, we may assume $y \geq x$ without loss of generality, using symmetry. Abbreviating $u_- := \max\{-u, 0\}$ and using $e^{-u} \geq 2u_-$, we get the claim (3.21):

(3.22)
$$\min\{x, y\} + e^{-x} + e^{-y} + \varepsilon|x - y|$$
$$\geq x + 2x_- + 2y_- + \delta|x - y|$$
$$= |x| - \delta x + 2y_- + \delta y$$
$$\geq (1 - \delta)|x| + \delta|y|$$
$$\geq \delta(|x| + |y|).$$

Combining (3.19), (3.20) and the definition (2.46) of $H_{\text{expI}}$, we conclude

(3.23)
$$H_{\text{middle},a,1/4} - H_{\text{expII}}$$
$$= (H_{\ln,a} + H_{\text{linear},a}) - (H_{\ln,1/2} + H_{\text{linear},1/2})$$
$$+ \left(H_{\ln,1/2} + H_{\text{linear},1/2} + H_{\text{tree}} + H_{\text{constraint}} - \frac{\Gamma}{4}\right) + H_{\text{expI}}$$

$$\geq \frac{\tilde{a}}{4}|\underline{X} + \overline{X} - \underline{X}' - \overline{X}' + 2W|$$

$$+ \frac{\tilde{a}}{8}[|\underline{X} - \overline{X}| + |\underline{X}' - \overline{X}'| + |\underline{X} + \overline{X} + \underline{X}' + \overline{X}' - 4Z|]$$

$$+ \frac{1}{4}[\min\{\underline{X}, \overline{X}\} + \min\{\underline{X}', \overline{X}'\} + e^{-\underline{X}} + e^{-\overline{X}} + e^{-\underline{X}'} + e^{-\overline{X}'}].$$

For the first summand in the last expression, we use the relation (2.41) between $W$ and $\Gamma$. Then, two applications of the fact (3.21) with $\varepsilon = \tilde{a}/2$, $\delta = \min\{\tilde{a}/2, 1/2\}$ yield the following lower bound for the terms (3.23):

$$\text{(3.23)} \geq \frac{\tilde{a}}{2}|\Gamma| + \frac{\tilde{a}}{8}|\underline{X} + \overline{X} + \underline{X}' + \overline{X}' - 4Z|$$

$$+ \frac{\delta}{4}[|\underline{X}| + |\overline{X}| + |\underline{X}'| + |\overline{X}'|]$$

(3.24)
$$\geq \frac{\delta}{8}|\Gamma| + \frac{\delta}{8}|\underline{X} + \overline{X} + \underline{X}' + \overline{X}' - 4Z|$$

$$+ \frac{\delta}{4}[|\underline{X}| + |\overline{X}| + |\underline{X}'| + |\overline{X}'|]$$

$$\geq \frac{\delta}{8}|\Gamma| + \frac{\delta}{2}|Z| + \frac{\delta}{8}[|\underline{X}| + |\overline{X}| + |\underline{X}'| + |\overline{X}'|]$$

$$\geq \frac{\delta}{8}[|\underline{X}| + |\overline{X}| + |Z| + |\Gamma| + |\underline{X}'| + |\overline{X}'|].$$

In the case $\eta = 1/4$, this proves the second inequality in the claim (3.17). By the symmetry property (2.49), this implies the second inequality in (3.17) in the case $\eta = -1/4$, too. For $-1/4 < \eta < 1/4$, $H_{\text{middle},a,\eta} - H_{\text{expII}}$ is a convex combination of $H_{\text{middle},a,1/4} - H_{\text{expII}}$ and $H_{\text{middle},a,-1/4} - H_{\text{expII}}$. Thus we get the second inequality in (3.17) in the general case, too. The first inequality in (3.17) follows immediately from $H_{\text{expII}} \geq 0$; recall its definition (2.47). □

LEMMA 3.3. *For all* $\omega = (\underline{X}, \overline{X}, \sigma, T) \in \Omega_{\text{cycle}}$, $\omega' = (\underline{X}', \overline{X}', \sigma', T') \in \Omega_{\text{cycle}}$, $\omega_{\text{rung}} = (Z, \Gamma) \in \Omega_{\text{rung}}$, *the map*

$$\text{(3.25)} \qquad [-1, 1] \ni \gamma \mapsto H_{\text{middle},a,0}(\omega | Z, \Gamma + \gamma \mathbb{1}_{\{\sigma \neq \sigma'\}} | \omega')$$

*is twice differentiable with the following bounds on its derivatives:*

(3.26)
$$\sup_{-1 \leq \gamma \leq 1} \left| \frac{\partial^j}{\partial \gamma^j} H_{\text{middle},a,0}(\omega | Z, \Gamma + \mathbb{1}_{\{\sigma \neq \sigma'\}} \gamma | \omega') \right|$$

$$\leq c_8 + H_{\text{expII}}(\omega | \omega_{\text{rung}} | \omega')$$

*holds for* $j = 1, 2$ *with some constant* $c_8 = c_8(a) > 0$.

PROOF. Let $j \in \{1,2\}$, and let $\Theta_\gamma$ denote the map $(\omega|\omega_{\mathrm{rung}}|\omega') \mapsto (\omega|Z, \Gamma + \mathbb{1}_{\{\sigma \neq \sigma'\}}\gamma|\omega')$.

By (2.41), $\Gamma = W + U - U'$. Hence, the terms contributing to the derivatives under consideration are those depending on $\Gamma$ and $W$. Consequently, using the explicit form of $H_{\mathrm{middle},a,0}$ from Definition 2.10, we see that

$$
\sup_{-1 \leq \gamma \leq 1} \left| \frac{\partial^j}{\partial \gamma^j}[H_{\mathrm{middle},a,0} \circ \Theta_\gamma] \right|
$$
$$
\leq c_9 + \sup_{-1 \leq \gamma \leq 1} \left| \frac{\partial^j}{\partial \gamma^j}[H_{\ln,a} \circ \Theta_\gamma] + \frac{\partial^j}{\partial \gamma^j}[H_{\mathrm{expII}} \circ \Theta_\gamma] \right| \tag{3.27}
$$

with a constant $c_9 > 0$. It is not hard to see that $\frac{\partial^j}{\partial W^j} \ln[e^{X+W/2} + e^{X'-W/2} + e^Z]$ is bounded for $j = 1, 2$. Hence,

$$
\sup_{-1 \leq \gamma \leq 1} \left| \frac{\partial^j}{\partial \gamma^j}[H_{\ln,a} \circ \Theta_\gamma] \right| \leq c_{10} \tag{3.28}
$$

with a constant $c_{10} > 0$. By the definition (2.47) of $H_{\mathrm{expII}}$, we have

$$
H_{\mathrm{expII}} \circ \Theta_\gamma = \tfrac{1}{2}[e^{(W \circ \Theta_\gamma)/2} - 2\sigma\sigma' + e^{-(W \circ \Theta_\gamma)/2}]e^{-Z} \tag{3.29}
$$

with $W \circ \Theta_\gamma = W + \mathbb{1}_{\{\sigma \neq \sigma'\}}\gamma$. Thus, the following hold:

$$
\sup_{-1 \leq \gamma \leq 1} \left| \frac{\partial}{\partial \gamma}[H_{\mathrm{expII}} \circ \Theta_\gamma] \right|
$$
$$
= \sup_{-1 \leq \gamma \leq 1} \left| \frac{1}{4}[e^{(W \circ \Theta_\gamma)/2} - e^{-(W \circ \Theta_\gamma)/2}]e^{-Z} \mathbb{1}_{\{\sigma \neq \sigma'\}} \right| \tag{3.30}
$$
$$
\leq H_{\mathrm{expII}},
$$

$$
\sup_{-1 \leq \gamma \leq 1} \left| \frac{\partial^2}{\partial \gamma^2}[H_{\mathrm{expII}} \circ \Theta_\gamma] \right|
$$
$$
= \sup_{-1 \leq \gamma \leq 1} \left| \frac{1}{8}[e^{(W \circ \Theta_\gamma)/2} + e^{-(W \circ \Theta_\gamma)/2}]e^{-Z} \mathbb{1}_{\{\sigma \neq \sigma'\}} \right| \tag{3.31}
$$
$$
\leq H_{\mathrm{expII}}.
$$

For the last inequalities in (3.30) and (3.31), we used that $|e^{(W+\gamma)/2} \pm e^{-(W+\gamma)/2}| \leq 2(e^{W/2} - 2\sigma\sigma' + e^{-W/2})$ for $-1 \leq \gamma \leq 1$ on the event $\{\sigma \neq \sigma'\}$, that is, for $-2\sigma\sigma' = 2$. The claim follows from the bounds (3.30) and (3.31) together with (3.27) and (3.28). $\square$

3.2. *Bounding $H_{\text{left},a}$ and $H_{\text{right},a}$.* The bounds for the boundary parts in the total Hamiltonian are obtained roughly similarly to the bound of the "middle" piece.

LEMMA 3.4. *Let $a = 3/4$. The following bounds hold:*

$$(3.32) \qquad H_{\text{left},3/4} - H_{\exp} \geq Z/4 \quad and \quad H_{\text{right},3/4} - H_{\exp} \geq Z/4.$$

PROOF. Let $T \in \{A, B, C, D\}$. First we estimate $H_{\text{left},3/4}$. For all $\alpha, \beta \in [0, 1]$, we have

$$(3.33) \qquad \begin{aligned} & (a + \tfrac{1}{2}) \ln[e^{\underline{X}} + e^Z] + a \ln[e^{\overline{X}} + e^Z] \\ & \geq \tfrac{5}{4} \max\{\underline{X}, Z\} + \tfrac{3}{4} \max\{\overline{X}, Z\} \\ & \geq \tfrac{5}{4}[\alpha \underline{X} + (1 - \alpha)Z] + \tfrac{3}{4}[\beta \overline{X} + (1 - \beta)Z]. \end{aligned}$$

We choose

$$(3.34) \qquad \alpha = \tfrac{1}{5}\mathbb{1}_{\{T=A\}} + \tfrac{3}{5}\mathbb{1}_{\{T=B\}} + \tfrac{2}{5}\mathbb{1}_{\{T=D\}},$$

$$(3.35) \qquad \beta = \tfrac{1}{3}\mathbb{1}_{\{T=A\}} + \mathbb{1}_{\{T=B\}} + \tfrac{2}{3}\mathbb{1}_{\{T=C\}}.$$

Substituting the bound (3.33) with this choice in Definition 2.11 of $H_{\text{left},3/4}$, an elementary but tedious calculation shows that the bound on the left-hand side in (3.32) is satisfied.

Similarly, we get a bound for $H_{\text{right},3/4}$: For all $\alpha, \beta \in [0, 1]$, we have

$$(3.36) \qquad \begin{aligned} & (a + \tfrac{1}{2})\{\ln[e^{\underline{X}} + e^Z] + \ln[e^{\overline{X}} + e^Z]\} \\ & \geq \tfrac{5}{4}[\max\{\underline{X}, Z\} + \max\{\overline{X}, Z\}] \\ & \geq \tfrac{5}{4}[\alpha \underline{X} + (1 - \alpha)Z + \beta \overline{X} + (1 - \beta)Z]. \end{aligned}$$

This time, we choose

$$(3.37) \qquad \alpha = \tfrac{4}{5}\mathbb{1}_{\{T=A\}} + \tfrac{2}{5}\mathbb{1}_{\{T=B\}} + \tfrac{1}{5}\mathbb{1}_{\{T=C\}} + \tfrac{3}{5}\mathbb{1}_{\{T=D\}},$$

$$(3.38) \qquad \beta = \alpha + \tfrac{2}{5}\mathbb{1}_{\{T=C\}} - \tfrac{2}{5}\mathbb{1}_{\{T=D\}}.$$

Substituting (3.36) in the definition of $H_{\text{right},3/4}$, the bound on the right-hand side in (3.32) is also satisfied. □

PROPOSITION 3.5. *For all $a > 3/4$, we have the following estimates:*

$$(3.39) \qquad H_{\text{left},a} \geq c_{11}(a)[|\underline{X}| + |\overline{X}| + |Z|],$$

$$(3.40) \qquad H_{\text{right},a} \geq c_{11}(a)[|\underline{X}| + |\overline{X}| + |Z|],$$

*with the constant $c_{11}(a) = \tfrac{1}{2}\min\{a - \tfrac{3}{4}, \tfrac{1}{6}\} > 0$.*

PROOF. Let $a > 3/4$, and let $\tilde{a} = a - 3/4$. Using the estimate $2\ln(e^x + e^y) - x - y \geq 2\max\{x, y\} - x - y = |x - y|$, we obtain

$$
\begin{aligned}
H_{\text{left},a} - H_{\text{left},3/4} &= \tilde{a}\{\ln[e^{\overline{X}} + e^Z] + \ln[e^{\underline{X}} + e^Z] - U - Z\} \\
&\geq \frac{\tilde{a}}{2}\{|\overline{X} - Z| + |\underline{X} - Z|\}.
\end{aligned}
\tag{3.41}
$$

Combining this with Lemma 3.4 and the bounds $H_{\exp} \geq e^{-Z}/2$ and $Z + 2e^{-Z} \geq |Z|$ yields

$$
\begin{aligned}
H_{\text{left},a} &\geq \frac{\tilde{a}}{2}\{|\overline{X} - Z| + |\underline{X} - Z|\} + \frac{1}{4}[Z + 2e^{-Z}] \\
&\geq c_{11}(a)\{|\overline{X} - Z| + |\underline{X} - Z| + 3|Z|\} \\
&\geq c_{11}(a)[|\underline{X}| + |\overline{X}| + |Z|].
\end{aligned}
\tag{3.42}
$$

This implies the claim (3.39). The estimate (3.40) for $H_{\text{right},a}$ follows with the same arguments; one just replaces "left" by "right." □

## 4. Statistical mechanics of the random environment.

4.1. *Transfer operators.* In this section we introduce a transfer operator $K_\eta$, and we show that it is a Hilbert–Schmidt operator. A Perron–Frobenius type argument yields all the spectral information about $K_\eta$ that we need.

DEFINITION 4.1. Let $d\omega_{\text{cycle}}$ denote the Lebesgue measure times the counting measure on $\Omega_{\text{cycle}}$. We define the following Hilbert space:

$$\mathcal{H} := L^2(\Omega_{\text{cycle}}, \mathcal{B}(\Omega_{\text{cycle}}), d\omega_{\text{cycle}}). \tag{4.1}$$

The scalar product in $\mathcal{H}$ is denoted by

$$\langle fg \rangle := \int_{\Omega_{\text{cycle}}} f(\omega_{\text{cycle}})g(\omega_{\text{cycle}})\, d\omega_{\text{cycle}}. \tag{4.2}$$

LEMMA 4.2. *Let $c_7(a)$ be as in Proposition 3.2, let $c_{11}(a)$ be as in Proposition 3.5 and let $c_{12}(a) := \min\{c_7(a), c_{11}(a)\}/2$. For $\omega_{\text{cycle}} = (\underline{X}, \overline{X}, \sigma, T) \in \Omega_{\text{cycle}}$ and $\omega_{\text{rung}} = (Z, \Gamma) \in \Omega_{\text{rung}}$, we define*

$$\|\omega_{\text{cycle}}\| := c_{12}(a)(|\underline{X}| + |\overline{X}|) \quad \text{and} \quad \|\omega_{\text{rung}}\| := c_{12}(a)(|Z| + |\Gamma|). \tag{4.3}$$

*Take a random variable $\Upsilon : \Omega_{\text{cycle}} \times \Omega_{\text{rung}} \times \Omega_{\text{cycle}} \to \mathbb{R}$ satisfying*

$$|\Upsilon(\omega|\omega_{\text{rung}}|\omega')| \leq c_{13}e^{\|\omega\| + \|\omega_{\text{rung}}\| + \|\omega'\| + H_{\text{expII}}(\omega|\omega_{\text{rung}}|\omega')} \tag{4.4}$$

for some constant $c_{13} > 0$. For $-1/4 \leq \eta \leq 1/4$, the function $k_\eta^\Upsilon : \Omega_{\text{cycle}} \times \Omega_{\text{cycle}} \to [0, \infty[$,

$$(4.5) \quad k_\eta^\Upsilon(\omega, \omega') := \int_{\Omega_{\text{rung}}} \Upsilon(\omega|\omega_{\text{rung}}|\omega') e^{-H_{\text{middle},a,\eta}(\omega|\omega_{\text{rung}}|\omega')} \, d\omega_{\text{rung}},$$

is well defined. The integral operator (acting from the right)

$$(4.6) \quad \begin{aligned} \tilde{K}_\eta^\Upsilon &: \mathcal{H} \to \mathcal{H}, \qquad f \mapsto f\tilde{K}_\eta^\Upsilon, \\ (f\tilde{K}_\eta^\Upsilon)(\omega') &= \int_{\Omega_{\text{cycle}}} f(\omega) k_\eta^\Upsilon(\omega, \omega') \, d\omega \end{aligned}$$

is a well-defined Hilbert–Schmidt operator. In the special case $\Upsilon = 1$, we write $k_\eta := k_\eta^1$ and $K_\eta := \tilde{K}_\eta^1$.

If we drop the summand $H_{\text{expII}}$ in the exponent of assumption (4.4), this assumption only gets stronger, because $H_{\text{expII}} \geq 0$ holds.

PROOF OF LEMMA 4.2. By Proposition 3.2, $H_{\text{middle},a,\eta} \geq 2(\|\omega\| + \|\omega_{\text{rung}}\| + \|\omega'\|) + H_{\text{expII}}$. Hence, using (4.4), the integrand in (4.5) is bounded by $c_{13} e^{-\|\omega\| - \|\omega_{\text{rung}}\| - \|\omega'\|}$. Integrating over $\omega_{\text{rung}}$, we get

$$(4.7) \quad |k_\eta^\Upsilon(\omega, \omega')| \leq c_{14} e^{-\|\omega\| - \|\omega'\|}$$

with a positive constant $c_{14}(a)$. In particular, $k_\eta^\Upsilon$ is well defined. Consequently, the following integral is finite:

$$(4.8) \quad \int_{\Omega_{\text{cycle}}} \int_{\Omega_{\text{cycle}}} k_\eta^\Upsilon(\omega, \omega')^2 \, d\omega \, d\omega' < \infty.$$

This shows that $\tilde{K}_\eta^\Upsilon : \mathcal{H} \to \mathcal{H}$ is a Hilbert–Schmidt operator. □

DEFINITION 4.3. When there is no risk of confusion, we use the following notation similar to the left and right operation of matrices: We denote the adjoint $K_\eta^* : \mathcal{H} \to \mathcal{H}$ of $K_\eta : \mathcal{H} \to \mathcal{H}$ by the same symbol $K_\eta$, but acting from the left:

$$(4.9) \quad \begin{aligned} K_\eta^* &: \mathcal{H} \to \mathcal{H}, \qquad g \mapsto K_\eta g, \\ (K_\eta g)(\omega) &= \int_{\Omega_{\text{cycle}}} k_\eta(\omega, \omega') g(\omega') \, d\omega'. \end{aligned}$$

In particular, $\langle fK_\eta g \rangle = \langle (fK_\eta)g \rangle = \langle f(K_\eta g) \rangle$.

LEMMA 4.4. *The spectral radius $\lambda_\eta > 0$ of $K_\eta$ is a simple eigenvalue of $K_\eta$ with unique (up to normalization) positive left and right eigenfunctions $v_\eta > 0$ and $v_\eta^* > 0$, that is, $v_\eta K_\eta = v_\eta \lambda_\eta$ and $K_\eta v_\eta^* = \lambda_\eta v_\eta^*$. Every other eigenvalue $\lambda$ of $K_\eta$ has modulus $|\lambda| < \lambda_\eta$.*

PROOF. This follows from Jentzsch's theorem; see, for example, Theorem 6.6 in [13] or Theorem 43.8 in [18]. We verify the hypotheses of Jentzsch's theorem: First, $K_\eta$ is compact by Lemma 4.2. Second, consider a subset $S \subset \Omega_{\text{cycle}}$ with the property that $S$ and $\Omega_{\text{cycle}} \setminus S$ have positive reference measure. We check that $\int_S \int_{\Omega_{\text{cycle}} \setminus S} k_\eta(\omega, \omega')\,d\omega\,d\omega' > 0$. Indeed, for $T \in \{A, B, C, D\}$, let $\Omega_T := \mathbb{R}^2 \times \{\pm 1\} \times \{T\} \subset \Omega_{\text{cycle}}$. Recall that $k_\eta(\omega, \omega') > 0$ iff $(\omega, \omega') \notin \Omega_A \times \Omega_B$. Note that $\Omega_C \cap S$ or $\Omega_C \setminus S$ has strictly positive reference measure, and $k_\eta(\omega, \omega') > 0$ for all $(\omega, \omega') \in (\Omega_C \cap S) \times (\Omega_{\text{cycle}} \setminus S)$ and for all $(\omega, \omega') \in S \times (\Omega_C \setminus S)$. Thus, the hypotheses of Jentzsch's theorem are fulfilled.

The same argument applies to $K_\eta^*$. Furthermore, it follows from the corollary to Proposition 5.1 on page 328 of [13] that $\lambda_\eta$ is an (algebraically) simple eigenvalue.

$\square$

We state an immediate consequence of this lemma:

COROLLARY 4.5. *Let $v_\eta$, $v_\eta^*$ be normalized such that $\langle v_\eta v_\eta^* \rangle = 1$. Let $M_\eta : \mathcal{H} \to \mathcal{H}$, $f \mapsto fM_\eta = \langle fv_\eta^* \rangle v_\eta$. Then as $n \to \infty$, $\lambda_\eta^{-n} K_\eta^n$ (acting from the right) converges to $M_\eta$ exponentially fast with respect to the operator norm $\|\cdot\|_{\mathcal{H} \to \mathcal{H}}$.*

DEFINITION 4.6. We define $g_{\text{left}}, g_{\text{right}} : \Omega_{\text{cycle}} \to \mathbb{R}$ by

$$g_{\text{left}}(\omega_{\text{cycle}}) := \int_{\mathbb{R}} e^{-H_{\text{left}}(Z, \omega_{\text{cycle}})}\,dZ, \tag{4.10}$$

$$g_{\text{right}}(\omega_{\text{cycle}}) := \int_{\mathbb{R}} e^{-H_{\text{right}}(\omega_{\text{cycle}}, Z)}\,dZ. \tag{4.11}$$

LEMMA 4.7. *One has $g_{\text{left}} \in \mathcal{H}$ and $g_{\text{right}} \in \mathcal{H}$.*

PROOF. By Proposition 3.5, $H_{\text{right},a} \geq \|\omega_{\text{cycle}}\| + c_{12}(a)|Z|$. Hence,

$$\|g_{\text{right}}\|_{\mathcal{H}}^2 = \int_{\Omega_{\text{cycle}}} \left[\int_{\mathbb{R}} e^{-H_{\text{right}}(\omega_{\text{cycle}}, Z)}\,dZ\right]^2 d\omega_{\text{cycle}} < \infty, \tag{4.12}$$

and $g_{\text{right}} \in \mathcal{H}$. Replacing in the above argument "right" by "left," we conclude that $g_{\text{left}} \in \mathcal{H}$.  $\square$

4.2. *A deformation of the Gibbs measure.* This section contains one of the central pieces of the whole proof of recurrence: We deform the Gibbs measure by changing $\eta = 1/4$ to the "more symmetric" value $\eta = 0$. In the language of statistical physics, the "spin chain" at the "physical" value $\eta =$ 1/4 is exposed to "external forces" of opposite directions at its two ends, causing the "symmetry breaking term" $\eta\Gamma_i$ in the Hamiltonian. The origin of this symmetry breaking term is the different scaling of the term $x_{\overline{0}}$ belonging to the starting point in the density (2.6). We compensate this external force artificially by setting $\eta = 0$, that is, by applying an external "counter-force," at least for the part of the spin chain between level 0 and level $j$. We define the corresponding "artificially deformed" Gibbs measure $\nu_{n,j}$:

LEMMA 4.8. *For $n \in \mathbb{N}$ and $j < n$, let*

$$\Sigma_j := \tfrac{1}{4} \sum_{i=1}^{j} \Gamma_i, \tag{4.13}$$

$$Z_{n,j} := E_{\mathbb{P}^{(n)}}[e^{-\Sigma_j}]. \tag{4.14}$$

*Then $Z_{n,j}$ is finite; thus the following probability measure is well defined:*

$$d\nu_{n,j} := \frac{e^{-\Sigma_j}}{Z_{n,j}} d\mathbb{P}^{(n)}. \tag{4.15}$$

PROOF. Using Definition 2.13, we know $Z_{n,j} Z^{(n)} = \int e^{-H^{(n)} - \Sigma_j} d\omega$. Note that

$$\begin{aligned} H^{(n)}(\omega) + \Sigma_j(\omega) \\ = H_{\text{left},a}(\omega_{\text{left}}|\omega_{\text{cycle},1}) \\ + \sum_{i=1}^{j} H_{\text{middle},a,0}(\omega_{\text{cycle},i}|\omega_{\text{rung},i}|\omega_{\text{cycle},i+1}) \\ + \sum_{i=j+1}^{n-1} H_{\text{middle},a,1/4}(\omega_{\text{cycle},i}|\omega_{\text{rung},i}|\omega_{\text{cycle},i+1}) \\ + H_{\text{right},a}(\omega_{\text{cycle},n}|\omega_{\text{right}}) \end{aligned} \tag{4.16}$$

holds; recall Definitions 2.10 and 2.12. As a consequence of the bounds (3.39), (3.40) and (3.17) for $H_{\text{left},a}$, $H_{\text{right},a}$ and $H_{\text{middle},a,\eta}$, the term $e^{-H^{(n)}(\omega) - \Sigma_j(\omega)}$ tends to 0 exponentially fast as at least one component of $\omega$ tends to $\pm\infty$. Thus $Z_{n,j}$ is finite, and we get

$$d\nu_{n,j} = \frac{e^{-H^{(n)}(\omega) - \Sigma_j(\omega)}}{Z_{n,j} Z^{(n)}} d\omega. \tag{4.17}$$

□

LEMMA 4.9. *Let $\Upsilon : \Omega_{\text{cycle}} \times \Omega_{\text{rung}} \times \Omega_{\text{cycle}} \to \mathbb{R}$ be a random variable such that the bound (4.4) holds. Then the expectation of $\Upsilon(\omega_{\text{cycle},i}|\omega_{\text{rung},i}|\omega_{\text{cycle},i+1})$*

with respect to $\nu_{n,j}$ exists and is uniformly bounded for $i, j, n$ with $i < n$ and $j < n$. It can be written as follows:

$$
\begin{aligned}
&E_{\nu_{n,j}}[\Upsilon(\omega_{\text{cycle},i}|\omega_{\text{rung},i}|\omega_{\text{cycle},i+1})] \\
(4.18) \quad &= \begin{cases} \dfrac{\langle g_{\text{left}} K_0^{i-1} \tilde{K}_0^{\Upsilon} K_0^{j-i} K_{1/4}^{n-j-1} g_{\text{right}}\rangle}{\langle g_{\text{left}} K_0^{j} K_{1/4}^{n-j-1} g_{\text{right}}\rangle}, & \text{for } i \le j < n, \\ \dfrac{\langle g_{\text{left}} K_0^{j} K_{1/4}^{i-j-1} \tilde{K}_{1/4}^{\Upsilon} K_{1/4}^{n-i-1} g_{\text{right}}\rangle}{\langle g_{\text{left}} K_0^{j} K_{1/4}^{n-j-1} g_{\text{right}}\rangle}, & \text{for } j < i < n. \end{cases}
\end{aligned}
$$

PROOF. We abbreviate $\hat{K}_\eta := \lambda_\eta^{-1} K_\eta$ and $\Upsilon_i := \Upsilon(\omega_{\text{cycle},i}|\omega_{\text{rung},i}|\omega_{\text{cycle},i+1})$. Using Fubini's theorem, (4.16), Lemma 4.2 and Definition 4.6, we calculate for $1 \le i \le j < n$:

$$
\begin{aligned}
E_{\nu_{n,j}}[\Upsilon_i] &= \frac{\int_{\Omega^{(n)}} \Upsilon_i e^{-H^{(n)}-\Sigma_j} d\omega}{\int_{\Omega^{(n)}} e^{-H^{(n)}-\Sigma_j} d\omega} = \frac{\langle g_{\text{left}} K_0^{i-1} \tilde{K}_0^{\Upsilon} K_0^{j-i} K_{1/4}^{n-j-1} g_{\text{right}}\rangle}{\langle g_{\text{left}} K_0^{j} K_{1/4}^{n-j-1} g_{\text{right}}\rangle} \\
(4.19) \\
&= \lambda_0^{-1} \frac{\langle g_{\text{left}} \hat{K}_0^{i-1} \tilde{K}_0^{\Upsilon} \hat{K}_0^{j-i} \hat{K}_{1/4}^{n-j-1} g_{\text{right}}\rangle}{\langle g_{\text{left}} \hat{K}_0^{j} \hat{K}_{1/4}^{n-j-1} g_{\text{right}}\rangle}.
\end{aligned}
$$

Note that the denominator in the last expression does not vanish. This shows that $E_{\nu_{n,j}}[\Upsilon_i]$ exists and is given by (4.18).

We claim that

$$
(4.20) \quad \inf_{j,m \in \mathbb{N}_0} \langle g_{\text{left}} \hat{K}_0^{j} \hat{K}_{1/4}^{m} g_{\text{right}}\rangle > 0.
$$

To prove this claim, we observe that for every fixed $j$ and $m$,

$$
(4.21) \quad \langle g_{\text{left}} \hat{K}_0^{j} \hat{K}_{1/4}^{m} g_{\text{right}}\rangle > 0
$$

holds, since this is a scalar product of positive functions. Furthermore, Corollary 4.5 implies

$$
(4.22) \quad \|g_{\text{left}} \hat{K}_0^{j} - \langle g_{\text{left}} v_0^* \rangle v_0\|_{\mathcal{H}} \le c_{15} e^{-c_{16} j}
$$

for some constants $c_{15}, c_{16} > 0$. As a consequence, we get for every fixed $m \ge 0$:

$$
(4.23) \quad \lim_{j \to \infty} \langle g_{\text{left}} \hat{K}_0^{j} \hat{K}_{1/4}^{m} g_{\text{right}}\rangle = \langle g_{\text{left}} v_0^* \rangle \langle v_0 \hat{K}_{1/4}^{m} g_{\text{right}}\rangle > 0;
$$

note that we have again scalar products of positive functions. Similarly, again by Corollary 4.5, we know $\hat{K}_{1/4}^{m} g_{\text{right}} \to v_{1/4}^* \langle v_{1/4} g_{\text{right}}\rangle$ in $\mathcal{H}$ (exponentially fast) as $m \to \infty$. Consequently, we get for every fixed $j \ge 0$:

$$
(4.24) \quad \lim_{m \to \infty} \langle g_{\text{left}} \hat{K}_0^{j} \hat{K}_{1/4}^{m} g_{\text{right}}\rangle = \langle g_{\text{left}} \hat{K}_0^{j} v_{1/4}^* \rangle \langle v_{1/4} g_{\text{right}}\rangle > 0.
$$

Finally, as $j$ and $m$ tend to $\infty$ simultaneously, we get

$$\lim_{\substack{j\to\infty \\ m\to\infty}} \langle g_{\text{left}} \hat{K}_0^j \hat{K}_{1/4}^m g_{\text{right}}\rangle = \langle g_{\text{left}} v_0^*\rangle \langle v_0 v_{1/4}^*\rangle \langle v_{1/4} g_{\text{right}}\rangle > 0. \quad (4.25)$$

Combining (4.21), (4.23), (4.24) and (4.25) yields our claim (4.20).

Furthermore, we claim that $\langle g_{\text{left}} \hat{K}_0^{i-1} \tilde{K}_0^\Upsilon \hat{K}_0^{j-i} \hat{K}_{1/4}^{n-j-1} g_{\text{right}}\rangle$ is uniformly bounded in $i, j, n$. Indeed, the sequence $(g_{\text{left}} \hat{K}_0^{i-1})_{i\in\mathbb{N}}$ is bounded in $\mathcal{H}$; it even converges. Similarly, $\hat{K}_0^{j-i} \hat{K}_{1/4}^{n-j-1} g_{\text{right}}$ is bounded in $\mathcal{H}$, too; recall that $\hat{K}_0$ and $\hat{K}_{1/4}$ have the leading simple eigenvalue 1. Thus the numerator on the right-hand side of (4.18) in the case $i \le j < n$ is uniformly bounded. Similarly, the expression $\langle g_{\text{left}} K_0^j K_{1/4}^{i-j-1} \tilde{K}_{1/4}^\Upsilon K_{1/4}^{n-i-1} g_{\text{right}}\rangle$, that is, the numerator on the right-hand side of (4.18) in the case $j < i < n$, is also uniformly bounded in $i, j, n$. Finally, the denominator is uniformly bounded away from 0. This completes the proof of the lemma. $\square$

We apply Lemma 4.9 to the random variables $e^{c_{17}Z_i}$, $e^{c_{17}\Gamma_i}$, $e^{c_{17}\underline{X}_i}$ and $e^{c_{17}\overline{X}_i}$ with some sufficiently small constant $c_{17}(a) > 0$. Then the exponential Chebyshev inequality yields immediately the following consequence:

COROLLARY 4.10. *There exist positive constants $c_{18}(a)$ and $c_{17}(a)$, such that for all $n \in \mathbb{N}$, $j < n$ and $M > 0$, one has*

$$\nu_{n,j}[|\Upsilon| \ge M] \le c_{18} e^{-c_{17}M} \quad (4.26)$$

*whenever $\Upsilon$ is any of the random variables $Z_i$, $\Gamma_i$, $\underline{X}_i$ or $\overline{X}_i$ with any admissible $i$. In particular, for $j = 0$ we have a bound for $\nu_{n,0} = \mathbb{P}^{(n)}$.*

PROOF OF THEOREM 2.3. We apply the transformations (2.17) and (2.58) to express probabilities with respect to $\tilde{\mathbb{Q}}^{(n)}$ in terms of $\mathbb{P}^{(n)}$. Using the expressions (2.35)–(2.40) for $Z_i$, $\Gamma_i$, $\underline{X}_i$ and $\overline{X}_i$, the theorem follows from Corollary 4.10. Note that $\ln|y_{i+1}/y_i| = -W_i/2 = -\Gamma_i/2 + \underline{X}_i/4 + \overline{X}_i/4 - \underline{X}_{i+1}/4 - \overline{X}_{i+1}/4$ and $\ln(z_i/y_i^2) = Z_i - W_i/2$ hold for $1 \le i \le n-1$, and that linear combinations of random variables with exponential tails have exponential tails as well. $\square$

The following lemma states the basic symmetry for $\eta = 0$: For $\eta = 0$, the "separation" $\Gamma$ between neighboring spins has expectation 0, at least "far away from the boundary of the ladder."

Recall the definition (4.6) of the integral operator $\tilde{K}_0^\Gamma$.

LEMMA 4.11. *We have $\langle v_0 \tilde{K}_0^\Gamma v_0^*\rangle = 0$.*

PROOF. Let $\eta \in [-1/4, 1/4]$. Recall Definition 2.10 of $H_{\text{middle},a,\eta}$ and its reflection symmetry property (2.49). Using the definition (4.5) of $k_\eta^\Upsilon$, this reflection symmetry yields $k_\eta^1(\omega, \omega') = k_{-\eta}^1(\omega'^{\leftrightarrow}, \omega^{\leftrightarrow})$ and $k_\eta^\Gamma(\omega, \omega') = -k_{-\eta}^\Gamma(\omega'^{\leftrightarrow}, \omega^{\leftrightarrow})$. Let $f \in \mathcal{H}$. We set $f^{\leftrightarrow}(\omega) := f(\omega^{\leftrightarrow})$. The symmetry properties of $k_\eta^1$ and $k_\eta^\Gamma$ imply

$$
(4.27) \quad \begin{aligned}
(fK_\eta)^{\leftrightarrow}(\omega') &= \int f(\omega) k_\eta^1(\omega, \omega'^{\leftrightarrow}) \, d\omega \\
&= \int f(\omega^{\leftrightarrow}) k_{-\eta}^1(\omega', \omega) \, d\omega \\
&= K_{-\eta}(f^{\leftrightarrow})(\omega')
\end{aligned}
$$

and

$$
(4.28) \quad \begin{aligned}
(f\tilde{K}_\eta^\Gamma)^{\leftrightarrow}(\omega') &= \int f(\omega) k_\eta^\Gamma(\omega, \omega'^{\leftrightarrow}) \, d\omega \\
&= -\int f(\omega^{\leftrightarrow}) k_{-\eta}^\Gamma(\omega', \omega) \, d\omega \\
&= -\tilde{K}_{-\eta}^\Gamma(f^{\leftrightarrow})(\omega').
\end{aligned}
$$

We apply (4.27) for the eigenfunction $v_\eta$: $\lambda_\eta v_\eta^{\leftrightarrow} = (v_\eta K_\eta)^{\leftrightarrow} = K_{-\eta}(v_\eta^{\leftrightarrow})$. Hence, using the uniqueness of the eigenfunctions up to normalization (Lemma 4.4), we conclude $v_{-\eta}^* = c_{19}(\eta) v_\eta^{\leftrightarrow}$ for some constant $c_{19}(\eta) > 0$. We calculate, using (4.28),

$$(4.29) \quad \langle v_0 \tilde{K}_0^\Gamma v_0^* \rangle = \langle (v_0^*)^{\leftrightarrow} (v_0 \tilde{K}_0^\Gamma)^{\leftrightarrow} \rangle = -\langle (v_0^*)^{\leftrightarrow} \tilde{K}_0^\Gamma v_0^{\leftrightarrow} \rangle = -\langle v_0 \tilde{K}_0^\Gamma v_0^* \rangle.$$

This yields the claim. □

The next lemma applies this reflection symmetry to finite ladders: In the symmetric case ($\eta = 0$), the separation $\Gamma_i$ between neighboring spins does not get too large, even in the presence of boundary terms, by an "approximate" reflection symmetry.

LEMMA 4.12. *There are positive constants $c_{20}(a)$, $c_{16}(a)$, $c_{21}(a)$ and $c_{22}(a)$, such that for all $i, j, n \in \mathbb{N}$ with $i \leq j < n$ the following bound holds (uniformly in $n$):*

$$(4.30) \quad |E_{\nu_{n,j}}[\Gamma_i]| \leq c_{20} e^{-c_{16} i} + c_{21} e^{-c_{22}(j-i)}.$$

PROOF. Recall that $\hat{K}_\eta = \lambda_\eta^{-1} K_\eta$. The representation (4.19) implies that

$$(4.31) \quad E_{\nu_{n,j}}[\Gamma_i] = \lambda_0^{-1} \frac{\langle g_{\text{left}} \hat{K}_0^{i-1} \tilde{K}_0^\Gamma \hat{K}_0^{j-i} \hat{K}_{1/4}^{n-j-1} g_{\text{right}} \rangle}{\langle g_{\text{left}} \hat{K}_0^j \hat{K}_{1/4}^{n-j-1} g_{\text{right}} \rangle}.$$

By (4.20), the denominator in the last expression is bounded away from zero. Thus it suffices to derive a bound for the numerator.

Recall that the sequence $(\hat{K}_{1/4}^m g_{\text{right}})_{m \in \mathbb{N}_0}$ is bounded. Using Corollary 4.5 once more, this implies for some constants $c_{22}(a), c_{23}(a), c_{24}(a) > 0$, uniformly in $m$:

$$
\begin{aligned}
\|\hat{K}_0^l \hat{K}_{1/4}^m g_{\text{right}} &- v_0^* \langle v_0 \hat{K}_{1/4}^m g_{\text{right}}\rangle\|_{\mathcal{H}} \\
&\leq \|(\hat{K}_0^*)^l - M_0^*\|_{\mathcal{H} \to \mathcal{H}} \|\hat{K}_{1/4}^m g_{\text{right}}\|_{\mathcal{H}} \\
&\leq c_{23} e^{-c_{22} l},
\end{aligned}
\tag{4.32}
$$

$$
\|v_0^* \langle v_0 \hat{K}_{1/4}^m g_{\text{right}}\rangle\|_{\mathcal{H}} \leq c_{24}. \tag{4.33}
$$

Using Lemma 4.11, we obtain

$$
\begin{aligned}
&|\langle g_{\text{left}} \hat{K}_0^{i-1} \tilde{K}_0^{\Gamma} \hat{K}_0^l \hat{K}_{1/4}^m g_{\text{right}}\rangle| \\
&= |\langle g_{\text{left}} \hat{K}_0^{i-1} \tilde{K}_0^{\Gamma} \hat{K}_0^l \hat{K}_{1/4}^m g_{\text{right}}\rangle - \langle g_{\text{left}} \hat{K}_0^{i-1} \tilde{K}_0^{\Gamma} v_0^*\rangle \langle v_0 \hat{K}_{1/4}^m g_{\text{right}}\rangle \\
&\quad + \langle g_{\text{left}} \hat{K}_0^{i-1} \tilde{K}_0^{\Gamma} v_0^*\rangle \langle v_0 \hat{K}_{1/4}^m g_{\text{right}}\rangle \\
&\quad - \langle g_{\text{left}} v_0^*\rangle \langle v_0 \hat{K}_0^{\Gamma} v_0^*\rangle \langle v_0 \hat{K}_{1/4}^m g_{\text{right}}\rangle| \\
&\leq \|g_{\text{left}} \hat{K}_0^{i-1}\|_{\mathcal{H}} \|\tilde{K}_0^{\Gamma}\|_{\mathcal{H} \to \mathcal{H}} \|\hat{K}_0^l \hat{K}_{1/4}^m g_{\text{right}} - v_0^* \langle v_0 \hat{K}_{1/4}^m g_{\text{right}}\rangle\|_{\mathcal{H}} \\
&\quad + \|g_{\text{left}} \hat{K}_0^{i-1} - \langle g_{\text{left}} v_0^*\rangle v_0\|_{\mathcal{H}} \|\tilde{K}_0^{\Gamma}\|_{\mathcal{H} \to \mathcal{H}} \|v_0^* \langle v_0 \hat{K}_{1/4}^m g_{\text{right}}\rangle\|_{\mathcal{H}}.
\end{aligned}
\tag{4.34}
$$

Combining this with (4.22), (4.32), (4.33) and using the facts that $(g_{\text{left}} \hat{K}_0^{i-1})_{i \in \mathbb{N}}$ is bounded in $\mathcal{H}$ and that $\tilde{K}_0^{\Gamma}$ is a bounded linear operator, we conclude:

$$
|\langle g_{\text{left}} \hat{K}_0^{i-1} \tilde{K}_0^{\Gamma} \hat{K}_0^l \hat{K}_{1/4}^m g_{\text{right}}\rangle| \leq c_{26} e^{-c_{22} l} + c_{25} e^{-c_{16}(i-1)}, \tag{4.35}
$$

uniformly in $m$, with some positive constants $c_{25}(a)$ and $c_{26}(a)$. Substituting the bounds (4.35) and (4.20) in (4.31) yields the claim (4.30). □

Intuitively speaking, the spin chain described by $\mathbb{P}^{(n)}$ is exposed to some external forces between level 0 and level $j$ that are missing in $\nu_{n,j}$. In order to estimate the effect of these external forces, we compare $\mathbb{P}^{(n)}$ with another "artificial" deformation $\tilde{\nu}_{n,j,\gamma}$ of $\nu_{n,j}$. This deformation is obtained by taking the image with respect to some "deformation map" $\Xi_{\gamma,j}$. Roughly speaking, one may view $\tilde{\nu}_{n,j,\gamma}$ as an artificial "caricature version" of $\mathbb{P}^{(n)}$, which is easier to compare with $\nu_{n,j}$ than $\mathbb{P}^{(n)}$ itself.

DEFINITION 4.13. For $\gamma \in \mathbb{R}$ and $j < n$, we define the bijection $\Xi_{\gamma,j} \colon \Omega^{(n)} \to \Omega^{(n)}$,

$$(4.36) \quad \Xi_{\gamma,j}(\omega) := (\omega_{\text{left}}, (\omega_{\text{cycle},i})_{i=1,\ldots,n}, (\tilde{\omega}_{\text{rung},i,j,\gamma})_{i=1,\ldots,n-1}, \omega_{\text{right}}),$$

where

$$(4.37) \quad \tilde{\omega}_{\text{rung},i,j,\gamma} := (Z_i, \Gamma_i + \gamma \mathbb{1}_{\{\sigma_i \neq \sigma_{i+1}, i \leq j\}}).$$

Thus only the components $\Gamma_i$ with $i \leq j$ get shifted by $\gamma$, but only in the case $\sigma_i \neq \sigma_{i+1}$; all other components are left unchanged. Let $\tilde{\nu}_{n,j,\gamma}$ denote the image measure of $\nu_{n,j}$ with respect to $\Xi_{\gamma,j}$.

By deforming the spin chain with $\Xi_{\gamma,j}$, a deformation of the separation $\Sigma_j$ between distant spins roughly proportional to $j$ is induced. This is shown in the following lemma.

LEMMA 4.14. *There are positive constants $c_{27}(a)$ and $c_{28}(a)$, such that for all $\gamma \in (0,1]$ and for all $n,j,i \in \mathbb{N}$ with $i \leq j < n$, we have the following bounds:*

$$(4.38) \quad |E_{\nu_{n,j}}[\Sigma_j]| \leq c_{27},$$

$$(4.39) \quad \nu_{n,j}[\sigma_i \neq \sigma_{i+1}] \geq c_{28},$$

$$(4.40) \quad E_{\tilde{\nu}_{n,j,\gamma}}[\Sigma_j] \geq c_{28} j\gamma - c_{27}.$$

PROOF. Using the definition (4.13) of $\Sigma_j$, the claim (4.38) follows immediately by summing the bound (4.30) over $i = 1, \ldots, j$. To prove (4.39), we consider the map $S_i \colon \Omega^{(n)} \to \Omega^{(n)}$,

$$(4.41) \quad S_i(\omega) := (\omega_{\text{left}}, (\tilde{\omega}_{\text{cycle},\iota,i})_{\iota=1,\ldots,n}, (\omega_{\text{rung},\iota})_{\iota=1,\ldots,n-1}, \omega_{\text{right}})$$

with

$$(4.42) \quad \tilde{\omega}_{\text{cycle},\iota,i} = (\underline{X}_\iota, \overline{X}_\iota, (\mathbb{1}_{\{\iota \leq i\}} - \mathbb{1}_{\{\iota > i\}})\sigma_\iota, T_\iota),$$

that is, just the sign components $\sigma_\iota$ in $\omega$ with $\iota > i$ get flipped; all other components are left unchanged. We calculate the Radon–Nikodym derivative of the image measure $S_i[\nu_{n,j}]$ with respect to $\nu_{n,j}$: The density of $\nu_{n,j}$ is given by (4.17). The part of the Hamiltonian depending on the signs $\sigma_\iota$ is $H_{\text{expII}}$; recall Definition 2.10, in particular (2.47). Thus, we have

$$(4.43) \quad \frac{dS_i[\nu_{n,j}]}{d\nu_{n,j}} = \frac{\exp[-H_{\text{expII}}(\sigma_i, W_i, Z_i, -\sigma_{i+1})]}{\exp[-H_{\text{expII}}(\sigma_i, W_i, Z_i, \sigma_{i+1})]}$$
$$= \exp[-2\sigma_i \sigma_{i+1} e^{-Z_i}].$$

ERRW ON A LADDER

We calculate:

$$\nu_{n,j}[\sigma_i \neq \sigma_{i+1}] = (S_i[\nu_{n,j}])[\sigma_i = \sigma_{i+1}]$$
$$(4.44) \qquad = E_{\nu_{n,j}}[\exp[-2\sigma_i\sigma_{i+1}e^{-Z_i}]\mathbb{1}_{\{\sigma_i=\sigma_{i+1}\}}],$$

hence

$$(4.45) \quad 2\nu_{n,j}[\sigma_i \neq \sigma_{i+1}] = E_{\nu_{n,j}}[\mathbb{1}_{\{\sigma_i\neq\sigma_{i+1}\}} + \exp[-2e^{-Z_i}]\mathbb{1}_{\{\sigma_i=\sigma_{i+1}\}}].$$

By Corollary 4.10, the distributions of the random variables $Z_i$ with respect to $\nu_{n,j}$ are tight, uniformly in $n$, $j$ and $i$. Hence, the claim (4.39) follows from (4.45).

Finally, in order to prove (4.40), we observe for $i \leq j$

$$(4.46) \qquad E_{\tilde{\nu}_{n,j,\gamma}}[\Gamma_i] = E_{\nu_{n,j}}[\Gamma_i + \gamma\mathbb{1}_{\{\sigma_i\neq\sigma_{i+1}\}}] \geq E_{\nu_{n,j}}[\Gamma_i] + c_{28}\gamma$$

by Definition 4.13 and (4.39). Hence, we get the claim (4.40) using (4.38):

$$(4.47) \qquad E_{\tilde{\nu}_{n,j,\gamma}}[\Sigma_j] \geq E_{\nu_{n,j}}[\Sigma_j] + c_{28}\gamma j \geq c_{28}\gamma j - c_{27}. \qquad \square$$

The following lemma considers deformations of the local Hamiltonians in the "separation variable" $\Gamma_i$. It is an ingredient for a deformation argument, below.

Recall the definition (4.37) of $\tilde{\omega}_{\text{rung},i,j,\gamma}$.

LEMMA 4.15. *For $i \leq j \leq n$, the map $f_{i,j,n}:[-1,1] \to \mathbb{R}$,*

$$(4.48) \qquad f_{i,j,n}(\gamma) := E_{\nu_{n,j}}[H_{\text{middle},a,0}(\omega_{\text{cycle},i}|\tilde{\omega}_{\text{rung},i,j,\gamma}|\omega_{\text{cycle},i+1})]$$

*is twice differentiable with the following bound on its derivative:*

$$(4.49) \qquad \sup_{-1\leq\gamma\leq 1} |f''_{i,j,n}(\gamma)| \leq 2c_{29}$$

*for some constant $c_{29} = c_{29}(a) > 0$.*

PROOF. By Lemma 3.3, the map $\gamma \mapsto H_{\text{middle},a,0}(\omega_{\text{cycle},i}|\tilde{\omega}_{\text{rung},i,j,\gamma}|\omega_{\text{cycle},i+1})$ is twice differentiable, and its derivatives satisfy

$$(4.50) \qquad \sup_{-1\leq\gamma\leq 1}\left|\frac{\partial^k}{\partial\gamma^k}H_{\text{middle},a,0}(\omega_{\text{cycle},i}|\tilde{\omega}_{\text{rung},i,j,\gamma}|\omega_{\text{cycle},i+1})\right|$$
$$\leq c_8 + H_{\text{expII}}(\omega_{\text{cycle},i}|\omega_{\text{rung},i}|\omega_{\text{cycle},i+1})$$

for $k = 1, 2$. Hence, by Lebesgue's dominated convergence theorem, it suffices to show that

$$(4.51) \qquad E_{\nu_{n,j}}[H_{\text{expII}}(\omega_{\text{cycle},i}|\omega_{\text{rung},i}|\omega_{\text{cycle},i+1})]$$

is uniformly bounded. Recall that $H_{\mathrm{expII}} \geq 0$ by its definition (2.47). Hence, the bound (4.4) holds for $\Upsilon = H_{\mathrm{expII}}$, and consequently, Lemma 4.9 implies that (4.51) is uniformly bounded in $i, j, n$. □

Next, we derive a bound for the relative entropy between $\nu_{n,j}$ and $\tilde{\nu}_{n,j}$, quadratic in $\gamma$, and linear in the length $j$.

LEMMA 4.16. *For all $\gamma \in [-1, 1]$ and for all $j, n \in \mathbb{N}$ with $j < n$, we have*

$$(4.52) \qquad E_{\tilde{\nu}_{n,j,\gamma}}\left[\ln \frac{d\tilde{\nu}_{n,j,\gamma}}{d\nu_{n,j}}\right] \leq c_{29} j \gamma^2$$

*with the constant $c_{29}(a)$ taken from Lemma 4.15.*

PROOF. From Definition 4.13, we know $\tilde{\nu}_{n,j,\gamma} = \Xi_{\gamma,j}[\nu_{n,j}]$ and $\Xi_{\gamma,j}^{-1} = \Xi_{-\gamma,j}$. This implies

$$(4.53) \qquad \frac{d\tilde{\nu}_{n,j,\gamma}}{d\nu_{n,j}} \circ \Xi_{\gamma,j} = \frac{d\nu_{n,j}}{d\tilde{\nu}_{n,j,-\gamma}}$$

and thus

$$(4.54) \qquad E_{\tilde{\nu}_{n,j,\gamma}}\left[\ln \frac{d\tilde{\nu}_{n,j,\gamma}}{d\nu_{n,j}}\right] = E_{\nu_{n,j}}\left[\ln \frac{d\nu_{n,j}}{d\tilde{\nu}_{n,j,-\gamma}}\right] =: h(\gamma).$$

Let us calculate the Radon–Nikodym derivative in the last expectation: Recall the representation (4.17) of $\nu_{n,j}$. Since the reference measure $d\omega$ on $\Omega^{(n)}$ is invariant with respect to the map $\Xi_{\gamma,j}$, we conclude

$$(4.55) \qquad d\tilde{\nu}_{n,j,-\gamma} = \frac{e^{-H^{(n)} \circ \Xi_{\gamma,j} - \Sigma_j \circ \Xi_{\gamma,j}}}{Z_{n,j} Z^{(n)}} \, d\omega.$$

Combining (4.17) and (4.55), and using Definition 4.13 and (4.16) for $H^{(n)} + \Sigma_j$, we get

$$\ln \frac{d\nu_{n,j}}{d\tilde{\nu}_{n,j,-\gamma}} = H^{(n)} \circ \Xi_{\gamma,j} + \Sigma_j \circ \Xi_{\gamma,j} - H^{(n)} - \Sigma_j$$

$$(4.56) \qquad = \sum_{i=1}^{j} [H_{\mathrm{middle},a,0}(\omega_{\mathrm{cycle},i}|\tilde{\omega}_{\mathrm{rung},i,j,\gamma}|\omega_{\mathrm{cycle},i+1})$$

$$- H_{\mathrm{middle},a,0}(\omega_{\mathrm{cycle},i}|\omega_{\mathrm{rung},i}|\omega_{\mathrm{cycle},i+1})].$$

Using Lemma 4.15, this implies that the relative entropy $h(\gamma)$ is twice differentiable with respect to $\gamma \in [-1, 1]$. Furthermore, $h(\gamma) \geq 0$ holds, and $h(0) = 0$ follows from $\tilde{\nu}_{n,j,0} = \nu_{n,j}$. Consequently $h(\gamma)$ reaches its minimum value 0 at $\gamma = 0$; thus $h'(0) = 0$ is valid. Using Taylor's formula, (4.56) and the bound (4.49), we conclude for some $\xi = \xi(\gamma) \in [0, 1]$

$$(4.57) \qquad h(\gamma) = \frac{\gamma^2}{2} h''(\xi\gamma) = \frac{\gamma^2}{2} \sum_{i=1}^{j} f''_{i,j,n}(\xi\gamma) \le c_{29} j \gamma^2.$$

Thus, the claim (4.52) follows from (4.54). □

The following proposition is the key ingredient to prove recurrence of reinforced random walks: It bounds exponential moments for the "separation variables" $\Sigma_j$. Here it becomes clear what the deformed measure $\tilde{\nu}_{n,j,\gamma}$ is good for: it just serves to bound the free energy difference between the $\nu_{n,j}$ and $\mathbb{P}^{(n)}$ via a variational principle for free energies.

PROPOSITION 4.17. *For some positive constants $c_{27}(a)$ and $c_1(a)$, we have the following bound for all $j, n \in \mathbb{N}$ with $j < n$:*

$$(4.58) \qquad E_{\mathbb{P}^{(n)}}[e^{-\Sigma_j}] \le e^{c_{27} - c_1 j}.$$

PROOF. From (4.14) and (4.15), we know

$$(4.59) \qquad \ln E_{\mathbb{P}^{(n)}}[e^{-\Sigma_j}] = \ln Z_{n,j} = -\ln \frac{d\nu_{n,j}}{d\mathbb{P}^{(n)}} - \Sigma_j.$$

We take $0 < \gamma \le 1$ sufficiently small (to be specified below). Using that the relative entropy between $\tilde{\nu}_{n,j,\gamma}$ and $\mathbb{P}^{(n)}$ is nonnegative, we conclude from (4.59):

$$
\begin{aligned}
\ln Z_{n,j} &\le E_{\tilde{\nu}_{n,j,\gamma}} \left[ \ln \frac{d\tilde{\nu}_{n,j,\gamma}}{d\mathbb{P}^{(n)}} \right] + \ln Z_{n,j} \\
&= E_{\tilde{\nu}_{n,j,\gamma}} \left[ \ln \frac{d\tilde{\nu}_{n,j,\gamma}}{d\mathbb{P}^{(n)}} - \ln \frac{d\nu_{n,j}}{d\mathbb{P}^{(n)}} - \Sigma_j \right] \\
&= E_{\tilde{\nu}_{n,j,\gamma}} \left[ \ln \frac{d\tilde{\nu}_{n,j,\gamma}}{d\nu_{n,j}} \right] - E_{\tilde{\nu}_{n,j,\gamma}}[\Sigma_j] \\
&\le c_{29} j \gamma^2 - c_{28} j \gamma + c_{27}.
\end{aligned}
$$
(4.60)

We used the bound (4.40) and Lemma 4.16 in the last step. Taking a fixed $\gamma > 0$ so small that $-c_1 := c_{29}\gamma^2 - c_{28}\gamma < 0$, the claim (4.58) follows. □

4.3. *Exponential decay of the random edge weights.* In this section we prove exponential bounds for the random environment.

PROOF OF THEOREM 1.2. Let $e \in E^{(n)}$ be an edge on level $i$ of the ladder, and let $c_1$ be as in Proposition 4.17. By Theorem 2.1 and the definition

(2.10) of $\tilde{\mathbb{Q}}^{(n)}$,

$$\begin{aligned}(4.61)\quad P^{(n)}\left[\lim_{t\to\infty}\frac{k_t(e)}{k_t(\{\underline{0},\overline{0}\})}\leq e^{-c_1 i}\right]&=\tilde{\mathbb{Q}}^{(n)}\left[\frac{x_e}{z_0}\leq e^{-c_1 i}\right]\\&=\mathbb{Q}^{(n)}\left[\frac{x_e}{z_0}\leq e^{-c_1 i}\right].\end{aligned}$$

Note that

$$(4.62)\qquad Y_i=-Z_0-\sum_{j=1}^{i-1}W_j=-Z_0-4\Sigma_{i-1}-U_i+U_1$$

by (2.30), (2.38) and (4.13). Using Lemma 2.13 and the transformation of variables (2.32)/(2.33), the right-hand side of (4.61) equals

$$(4.63)\qquad (4.61)=\mathbb{P}^{(n)}[e^{-4\Sigma_{i-1}+\Upsilon_i}\leq e^{-c_1 i}],$$

where $\Upsilon_i=\underline{X}_i-Z_0-U_i+U_1$ if $x_e=\underline{x}_i$, $\Upsilon_i=\overline{X}_i-Z_0-U_i+U_1$ if $x_e=\overline{x}_i$, $\Upsilon_i=Z_n-Z_0+U_1-U_n$ if $x_e=z_n$ and $i=n$, and $\Upsilon_i=Z_i-W_i/2-Z_0-U_i+U_1=Z_i-Z_0+U_1-(\Gamma_i+U_i+U_{i+1})/2$ if $x_e=z_i$ with $i\in\{1,2,\ldots,n-1\}$; the last case follows from (2.39) and (2.29). Recall that $U_i=(\underline{X}_i+\overline{X}_i)/2$ by (2.28). We use (4.58) and the exponential Chebyshev inequality to obtain

$$(4.64)\qquad \mathbb{P}^{(n)}[-4\Sigma_{i-1}\geq -2c_1 i]\leq e^{c_{27}-c_1(i-2)/2}.$$

Applying Corollary 4.10 with $j=0$ yields

$$(4.65)\qquad \mathbb{P}^{(n)}[\Upsilon_i\geq c_1 i]\leq c_{30}e^{-c_{31}i}$$

with some positive constants $c_{30}(a),c_{31}(a)>0$. Combining the bounds (4.64) and (4.65) with (4.63), we obtain the claim (1.6):

$$(4.66)\qquad (4.61)\geq 1-e^{c_{27}-c_1(i-2)/2}-c_{30}e^{-c_{31}i}.\qquad\square$$

**5. Recurrence.** In this section we prove Theorem 1.1. Recall that we call the reinforced random walk recurrent if almost all paths visit all vertices infinitely often.

LEMMA 5.1. *On any graph $G$, the edge-reinforced random walk is recurrent if almost all paths return to the starting point infinitely often.*

PROOF. A Borel–Cantelli argument shows that for any $u,v\in V$ with $\{u,v\}\in E$ the following holds:

$$(5.1)\quad P_{v_0,a}^G[u\text{ is visited infinitely often and }v\text{ is visited at most finitely often}]=0.$$

(Details can be found in the proof of Corollary 3.1 in [12].) Hence, the claim follows by induction. $\square$

LEMMA 5.2. *For all $a > 3/4$, the edge-reinforced random walk on $\mathbb{N}_0 \times \{1,2\}$ starting in $\overline{0}$ with all initial weights equal to $a$ is recurrent.*

PROOF. For $x \in \mathbb{R}_+^E$, let $Q_x$ denote the distribution of the nonreinforced random walk on $\mathbb{N}_0 \times \{1,2\}$ starting at $\overline{0}$ which jumps from $u$ to $v$ with probability proportional to the weight $x_{\{u,v\}}$ of the edge $\{u,v\}$. Let $A_k^{(n)}$ denote the event that the random walker returns $\geq k$ times to $\overline{0}$ before hitting the set $\{\underline{n}, \overline{n}\}$. We set $A^{(n)} := A_1^{(n)}$.

Our proof uses the connection between random walks and electric networks. Consider the graph $\mathbb{N}_0 \times \{1,2\}$ as an electric network, where edge $e$ has the conductance $x_e$ or equivalently the resistance $x_e^{-1}$. Let $R^{(n)}(x)$ denote the resistance between $\overline{0}$ and the set $\{\underline{n}, \overline{n}\}$ of the finite ladder $G^{(n)}$. Then $C^{(n)}(x) = 1/R^{(n)}(x)$ is the effective conductance. Recall that $x_{\overline{0}} = z_0 + \overline{x}_1 = 1 + \overline{x}_1 \geq 1$ $\mathbb{Q}^{(n)}$-a.s. The escape probability $Q_x[(A^{(n)})^c]$ can be expressed as follows (see, e.g., Section 1.3.4 of [5]):

$$(5.2) \qquad Q_x[(A^{(n)})^c] = \frac{C^{(n)}(x)}{x_{\overline{0}}} \leq C^{(n)}(x) = \frac{1}{R^{(n)}(x)} \qquad \mathbb{Q}^{(n)}\text{-a.s.}$$

By Rayleigh's monotonicity law (Section 1.4 of [5]), the effective resistance of an electric network decreases if some of the individual resistances decrease. In particular, the network obtained from the ladder $G^{(n)}$ by shorting together each of the vertices $\{\underline{i}, \overline{i}\}$ has an effective resistance $\tilde{R}^{(n)}(x) \leq R^{(n)}(x)$. We calculate

$$(5.3) \qquad \tilde{R}^{(n)}(x) = \sum_{i=1}^n \frac{1}{\underline{x}_i + \overline{x}_i} \geq \frac{1}{\underline{x}_n + \overline{x}_n}.$$

Denoting by $b \wedge c$ the minimum of $b$ and $c$, we conclude

$$(5.4) \qquad \begin{aligned} Q_x[(A^{(n)})^c] &\leq \frac{1}{R^{(n)}(x)} \wedge 1 \\ &\leq \frac{1}{\tilde{R}^{(n)}(x)} \wedge 1 \leq (\underline{x}_n + \overline{x}_n) \wedge 1 \qquad \mathbb{Q}^{(n)}\text{-a.s.} \end{aligned}$$

By the strong Markov property, $Q_x[A_k^{(n)}] = (Q_x[A^{(n)}])^k$. Using (2.18), we obtain

$$(5.5) \qquad \begin{aligned} P^{(n)}[A_k^{(n)}] &= \int_{\Lambda^{(n)}} Q_x[A_k^{(n)}] \mathbb{Q}^{(n)}(dx\,dy\,dT) \\ &\geq \int_{\Lambda^{(n)}} (Q_x[A^{(n)}])^k \mathbb{Q}^{(n)}(dx\,dy\,dT) \\ &\geq \int_{\Lambda^{(n)}} [1 - [(\underline{x}_n + \overline{x}_n) \wedge 1]]^k \mathbb{Q}^{(n)}(dx\,dy\,dT). \end{aligned}$$

Let $B_n := \{\ln \underline{x}_n \leq -c_1 n \text{ and } \ln \overline{x}_n \leq -c_1 n\}$. Theorem 2.1 and Lemma 2.6 imply that the distribution of the limiting vector $(\lim_{t \to \infty} k_t(e)/k_t(\{\underline{0}, \overline{0}\}))_e$ under $P^{(n)}$ equals $\mathbb{Q}^{(n)}$ because $z_0 = 1$ $\mathbb{Q}^{(n)}$-a.s. Hence, we conclude from Theorem 1.2 that

$$P^{(n)}[A_k^{(n)}] \geq \int_{B_n} [1 - [(\underline{x}_n + \overline{x}_n) \wedge 1]]^k \mathbb{Q}^{(n)}(dx\, dy\, dT)$$

(5.6)
$$\geq \mathbb{Q}^{(n)}[B_n](1 - 2e^{-c_1 n} \wedge 1)^k$$
$$\geq (1 - 2c_2 e^{-c_3 n})(1 - 2e^{-c_1 n} \wedge 1)^k.$$

The probability that the reinforced random walk starting in $\overline{0}$ does not return to $\overline{0}$ before hitting $\{\underline{n}, \overline{n}\}$ is the same for the finite ladder $G^{(n)}$ and for $\mathbb{N}_0 \times \{1, 2\}$. Consequently, if we denote by $P^+$ the distribution of the edge-reinforced random walk on $\mathbb{N}_0 \times \{1, 2\}$ starting at $\overline{0}$ with all initial edge weights equal to $a$, then $P^+[A_k^{(n)}] = P^{(n)}[A_k^{(n)}]$, and we conclude

(5.7)
$$P^+[\text{return } \geq k \text{ times to } \overline{0}] = P^+\left[\bigcup_{n=1}^{\infty} A_k^{(n)}\right]$$
$$= \lim_{n \to \infty} P^{(n)}[A_k^{(n)}] = 1.$$

Since $k$ is arbitrary, we conclude that the reinforced random walk on the half-ladder returns to $\overline{0}$ infinitely often with probability 1. The claim follows from Lemma 5.1. □

PROOF OF THEOREM 1.1. Because of (5.1) and Lemma 5.1, it is enough to show that the set $\{\underline{0}, \overline{0}\}$ is visited infinitely often.

For $v_0 \in V$ and $x = (x_e)_{e \in E} \in \mathbb{R}_+^E$, we denote by $P_{v_0, x}$ the distribution of the edge-reinforced random walk on $\mathbb{Z} \times \{1, 2\}$ starting in $v_0$ with initial edge weights given by $x$. If $x_e = a$ for all $e$, we write simply $P$. The distributions of the edge-reinforced random walk on $\mathbb{N}_0 \times \{1, 2\}$ and $-\mathbb{N}_0 \times \{1, 2\}$ are denoted by $P^+$ and $P^-$, respectively.

Let $\tau_0 := 0$, and for $i \geq 0$, let $\tau_i$ be the $i$th return time to the set $\{\underline{0}, \overline{0}\}$. We denote by $\rho_i$ the hitting time of the set $\{\underline{-1}, \overline{-1}, \underline{1}, \overline{1}\}$ after time $\tau_i$. By (5.1), $\rho_i < \infty$ $P$-a.s. on the set $\{\tau_i < \infty\}$.

We prove by induction on $i \in \mathbb{N}_0$ that $\tau_i < \infty$ holds $P$-a.s. This is clearly true for $i = 0$. For the induction step from $i$ to $i + 1$, assume that the claim holds for $i$. Then, $\rho_i < \infty$ $P$-a.s., and we obtain

(5.8)
$$P[\tau_{i+1} < \infty] = P[\tau_{i+1} < \rho_i]$$
$$+ P[\tau_{i+1} > \rho_i] \cdot P[\tau_{i+1} < \infty \mid \tau_{i+1} > \rho_i]$$

$$= P[\tau_{i+1} < \rho_i]$$
$$+ P[\tau_{i+1} > \rho_i] \cdot \int P_{X_{\rho_i}, w_{\rho_i}}[\text{visit } \{\underline{0}, \overline{0}\}] \, dP.$$

Let us consider a fixed realization of $X_{\rho_i}$ and $w_{\rho_i}$. Since $P_{X_{\rho_i}, w_{\rho_i}}[\text{visit } \{\underline{0}, \overline{0}\}]$ depends only on the edge weights on one of the two half-ladders $\mathbb{N}_0 \times \{1, 2\}$ or $-\mathbb{N}_0 \times \{1, 2\}$, this probability agrees with the corresponding one for the reinforced random walk on the half-ladder. Hence, for some $* \in \{-, +\}$ and some finite path $(v_0 = \overline{0}, v_1, \ldots, v_t = X_{\rho_i})$ we have

(5.9) $\quad P_{X_{\rho_i}, w_{\rho_i}}[\text{visit } \{\underline{0}, \overline{0}\}]$
$\quad = P^*[\text{visit}\{\underline{0}, \overline{0}\} \text{ after time } t | X_s = v_s \text{ for all } s \leq t].$

By Lemma 5.2, the reinforced random walk on $\mathbb{N}_0 \times \{1, 2\}$ is recurrent. By symmetry, the same is true for the reinforced random walk on $-\mathbb{N}_0 \times \{1, 2\}$. In particular, for these two processes, almost all paths return to $\{\underline{0}, \overline{0}\}$. Hence, the conditional probability on the right-hand side of (5.9) equals 1. It follows from (5.8), that $P[\tau_{i+1} < \infty] = 1$, that is, the reinforced random walk on $\mathbb{Z} \times \{1, 2\}$ returns to $\{\underline{0}, \overline{0}\}$ at least $i+1$ times. By the induction principle, there are $P$-a.s. infinitely many visits to the set $\{\underline{0}, \overline{0}\}$, and the claim follows. $\square$

**Acknowledgment.** S. Rolles would like to thank Mike Keane for getting her interested in the problem and sharing his enthusiasm about reinforced random walks.

MATHEMATICAL INSTITUTE  
UNIVERSITY OF MUNICH  
THERESIENSTR. 39  
D-80333 MUNICH  
GERMANY  
E-MAIL: merkl@mathematik.uni-muenchen.de  
URL: www.mathematik.uni-muenchen.de/˜merkl  

FACULTY OF MATHEMATICS  
UNIVERSITY OF BIELEFELD  
POSTFACH 100131  
D-33501 BIELEFELD  
GERMANY  
E-MAIL: srolles@math.uni-bielefeld.de  
URL: www.math.uni-bielefeld.de/˜srolles